\documentclass[11pt,reqno]{amsart}

\usepackage{amsfonts}
\usepackage{eurosym}
\usepackage{amssymb}
\usepackage{amsthm}
\usepackage{amsmath}
\usepackage{bm}
\usepackage{cite}
\usepackage{mathrsfs}
\usepackage{xcolor}
\usepackage[OT1]{fontenc}
\usepackage[left=2cm, right=2cm, top=3cm]{geometry}
\usepackage{hyperref}
\hypersetup{colorlinks=true, linkcolor=blue, citecolor=red}

\usepackage{times}

\usepackage[shortlabels]{enumitem}
\usepackage{enumitem}
\newcommand{\subscript}[2]{$#1 _ #2$}
\numberwithin{equation}{section}
\usepackage[square,numbers,sectionbib]{natbib}
\usepackage{varioref}
\usepackage{listings}
\usepackage{color}
\definecolor{dkgreen}{rgb}{0,0.6,0}
\definecolor{gray}{rgb}{0.5,0.5,0.5}
\definecolor{mauve}{rgb}{0.58,0,0.82}

\allowdisplaybreaks
\newtheorem{theorem}{Theorem}[section]

\newtheorem{definition}[theorem]{Definition}
\newtheorem{lemma}[theorem]{Lemma}

\theoremstyle{remark}
\newtheorem{remark}[theorem]{Remark}

\def\ds{d\sigma}

\def\V{{\mathbf{V}}}
\def\H{\mathbf{H}}
\def\W{\mathbf{W}}
\def\A{\mathbf{A}}
\def\L{\mathbf{L}}

\def\g{\textbf{\textit{g}}}
\def\ddt{\frac{\d}{\d t}}
\def\C{\mathcal{C}}

\def\n{\mathbf{n}}
\def\E{\mathcal E}

\newcommand{\numberset}{\mathbb}
\newcommand{\N}{\numberset{N}}
\newcommand{\R}{\numberset{R}}
\def\nablag{\nabla_\Gamma}
\def\divg{\mathrm{div}_\Gamma}
\def\J{\mathbf{J}}
\def\cd{\cdot}

\def\vn{v_\n}
\def\intg{\int_{\Gamma(t)}}
\def\ints#1{\int_{#1(t)}}
\def\Vn{v_\n}
\def\VVn{\V_\n}

\def\dt{\partial^\circ_t}
\def\P{\mathbf{P}}
\def\v{\mathbf{v}}
\def\ig{\intg}
\def\u{\mathbf{u}}
\def\vphi{\varphi}

\def\tphi{\widetilde{\phi}}
\def\eps{\epsilon}

\def\A{\mathbf{A}}

\def\phit_{\phi_t}

\def\tphimt{\tphi_{-t}}

\def\dr{(\rho_1-\rho_2)}
\def\ro{\rho}
\def\tphimt{\tphi_{-t}}

\def\tT{\widetilde{T}}
\def\H{\mathbf{H}}

\def\w{\mathbf{w}}

\def\norm#1{\left\Vert#1\right\Vert }
\def\norma#1{\left\vert#1\right\vert }
\def\trho{\rho}
\def\f{\mathbf f}

\def\ddt{\frac{d}{dt}}

\def\gz{\Gamma_0}
\def\LLQ(#1,#2){L^{#1}(0,\widetilde{T};\L^{#2}(\Gamma_0))}

\def\dts{\partial_t^*}

\def\gam{\Gamma}
\def\gt{(\gam(t))}
\def\NN{\mathcal{N}}
\def\Linfg{\L^\infty\gt}
\def\wu{\widehat{\u}}
\def\na{\nablag}
\def\nus{\frac{\nu_*}{16}\norm{\E_S(\V)}^2}
\def\non{\nonumber}
\def\Phitn{\Phi_t^n}
\def\z{\mathbf z}
\def\X{\mathbf X}
\everymath{\displaystyle}

\def \no#1#2#3 {{\bf #1} (#3), #2.}
\def \eds#1#2#3 {#1, #2, #3.}

\author[H.~Abels]{Helmut Abels}\author[H.~Garcke]{Harald Garcke}
\address[H.~Abels, H.~Garcke]{Fakult\"{a}t f\"{u}r Mathematik\\
Universit\"{a}t Regensburg\\
93040 Regensburg, Germany}
\email[H.~Abels]{helmut.abels@ur.de}
\email[H.~Garcke]{harald.garcke@ur.de}
\author[A.~Poiatti]{Andrea Poiatti}
\address[A.~Poiatti]{Faculty of Mathematics, University of Vienna, Oskar-Morgenstern-Platz
1, A-1090 Vienna, Austria}
\email[A.~Poiatti]{andrea.poiatti@univie.ac.at}

\begin{document}

\title[Diffuse Interface Model for Two-Phase Flows on Evolving Surfaces]{
Diffuse Interface Model for Two-Phase Flows on Evolving Surfaces with Different Densities: Global Well-Posedness
}

\keywords{Diffuse interface models, PDEs on surfaces, Navier-Stokes/Cahn-Hilliard equation}

\subjclass[2010]{35Q30,76D05,  35D35, 35Q35, 76T06}


\date{\today}
\maketitle

\begin{abstract}
  We show existence and uniqueness of strong solutions to a Navier-Stokes/Cahn-Hilliard type system on a given two-dimensional evolving surface in the case of different densities and a singular (logarithmic) potential. The system describes a diffuse interface model for a two-phase flow of viscous incompressible fluids on an evolving surface. We also establish the  validity of the instantaneous strict separation property from the pure phases. To show these results we use our previous achievements on local well-posedness together with suitable novel regularity results for the convective Cahn-Hilliard equation. The latter allows to obtain higher-order energy estimates to extend the local solution globally in time. To this aim the time evolution of energy type quantities has to be calculated and estimated carefully. 
\end{abstract}
\section{Introduction}

Two-phase flows and phase separation phenomena play an important role in the sciences and engineering applications. In the present contribution we continue our study in \cite{AGP1} of a diffuse interface model for the two-phase of incompressible fluids on an evolving interface $\mathcal{S}:=\bigcup_{t\in[0,T]}\Gamma(t)\times \{t\}$, $T>0$, for different densities in the case that the evolving interface is given. The model leads to the following Navier-Stokes/Cahn-Hilliard system for $(\v,\pi,\vphi,\mu)$  on $S$: 
 \begin{align}
     \begin{cases}
     \rho \P\dt \v +((\rho\v+\J_\rho)\cdot\nabla_\Gamma)\v+\rho \Vn\mathbf{H}\v+\Vn\H\J_\rho-2\P\divg(\nu(\vphi)\E_S(\v))+\nablag \pi\\
     =-\P\divg(\nablag\vphi\otimes\nablag\vphi)+2\P\divg(\nu(\vphi)\Vn\mathbf{H})+\frac\rho2\nablag(\Vn)^2,\\
         \divg\v=-H\Vn,
         \\\dt\vphi-\divg(m(\vphi)\nablag\mu)+\nablag\vphi\cdot \v=0,\\
         \mu=-\Delta_\Gamma\vphi+\Psi'(\vphi),
     \end{cases}
     \label{mainp}
 \end{align}
 with $\J_\rho=-\frac{\widetilde{\rho}_1-\widetilde{\rho}_2}{2}\nablag\mu$. Here   $\VVn=\Vn\n$ is the normal velocity of $\Gamma$ satisfying the compatibility assumption
 \begin{equation*}
  \intg H\Vn\ds=0 \qquad \text{for all }t\in [0,T],  
 \end{equation*}
   $H$ and $\mathbf{H}$ at time $t$ are the mean curvature and the Weingarten map of $\Gamma(t)\subseteq \R^3$ and $\P(x,t)$ is the orthogonal projection onto $T_x\Gamma(t)$ for $(x,t)\in \mathcal{S}$. These quantities are determined by the given evolving surface $\mathcal{S}$. Moreover, $\v\colon \mathcal{S}\to \R^3$ with $\v(x,t)\in T_{x}\Gamma(t)$ for all $(x,t)\in \mathcal{S}$ is the tangential (mean) velocity of the fluid mixture,  $\pi\colon \mathcal{S} \to \R$ is its pressure, and $\E_S(\v)$ is the surface
rate-of-strain tensor, defined as the symmetrized tangential velocity gradient of $\v$. Furthermore, $\varphi \colon \mathcal{S} \to [-1,1]$ is the volume fraction difference of the fluids, $\mu\colon \mathcal{S} \to \R$ is the chemical potential and $\Psi \colon [-1,1]\to\R$ a homogeneous free energy density  of the mixture. Also,   $\rho=\rho(\varphi):=\tfrac{\tilde\rho_1+\tilde\rho_2}2+ \tfrac{\tilde\rho_1-\tilde\rho_2}2\varphi$ is the mass density of the mixture, where $\tilde\rho_j>0$ is the specific density of fluid $j=1,2$. For simplicity we will set the mobility $m(\vphi)\equiv 1$. More details on the notation are given in Section~\ref{sec:Assump} below. We refer to \cite{AGP1} for the derivation of this model and references on related models and their analysis. 
We note that, as obtained in the derivation of the model in \cite{AGP1}, the density $\rho$ satisfies 
\begin{align}
\label{ma}
\dt\rho+\divg(\J_\rho)+\divg(\rho\u)=0.
\end{align}
In \cite{AGP1} existence and uniqueness of strong solutions locally in time for sufficiently regular initial data was shown. It is the purpose of this contribution to extend this result and show existence of strong solution globally in time in the case of singular potential $\Psi$, as well as to establish the validity of the instantaneous strict separation property of $\vphi$ from the pure phases. If the Navier-Stokes/Cahn-Hilliard system \eqref{mainp} is considered in two-dimensional Euclidean space instead of a two-dimensional evolving surface $\Gamma(t)\subseteq \mathbb{R}^3$, the corresponding result was shown by Giorgini~\cite{GiorginiWellposednessTwoDim} in the case of periodic boundary conditions and by Abels, Garcke, and Giorgini in \cite{AGGio} in the case of a bounded sufficiently smooth domain $\Omega \subseteq \R^2$. In the latter reference also results on the case of three space dimensions and further references for Navier-Stokes/Cahn-Hilliard systems in Euclidean space can be found. In bounded domains we also refer to Abels, Garcke, and Poiatti \cite{AGP} for the extension of the analysis to multi-phase flows and to \cite{GGGP} for its nonlocal counterpart. We note that in the case of same densities, i.e., $\rho\equiv const.$, global well-posedness of weak solutions was shown in Elliott and Sales \cite{ES} by combining some techniques developed by Caetano and Elliott~\cite{DE} and Olshanskii, Reusken, and Zhiliakov~\cite{V2}. In the latter contribution global existence and uniqueness of weak solutions for a (tangential) Navier-Stokes system on an evolving given surface is shown (we refer to \cite{BRS} for a review on the modeling of tangential Navier-Stokes equations on evolving surfaces).

In order to prove our main result on the existence  of global strong solutions we will use similar arguments as in \cite{AGGio}, but in the present situation all the results and arguments have to be generalized to the more delicate case of an evolving surface $\Gamma(t)$. To this end we also use some techniques and identities as by Caetano, Elliott, Grasselli, and Poiatti \cite{CEGP} for the Cahn-Hilliard equation with a singular potential on an evolving surface. The structure of the contribution is as follows. In Section~\ref{sec:Assump} we summarize some basic notation, assumptions and preliminary results. The main result is stated in Section~\ref{sec:main}. The first important step for the proof of the main result is to generalize the regularity results for the convective Cahn-Hilliard equation shown in \cite{AGGio} to the case of an evolving surface, which is done in Section~\ref{sec:ConvCH}. Finally, in Section~\ref{sec:proof} the main result is proven.

\section{Assumptions and Functional Setting}\label{sec:Assump}
\subsection{Flow Map and Surface Piola Transform}
\label{regflowmap1}
We use the same notation and assumptions as in \cite{AGP1}. In particular we assume that $\mathcal{S}= \bigcup_{0\leq t\leq T}\{t\}\times\Gamma(t)$ is a $C^6$-manifold. Moreover, we have a normal flow map
$\Phi^n_{(\cd)}\in C^4(\mathcal{S}_0,\mathcal{S})$, where $\mathcal{S}_0= [0,T_0]\times \Gamma_0$,  which solves
\begin{align*}
  &\ddt \Phi_t^n(\z)=\VVn(t,\Phi_t^n(\z)),\\&
\Phi_0^n(\z)=\z,
\end{align*}
for all $\z\in \Gamma_0$, $t\in[0,T]$, where ${\n}\in C^5([0,T]\times \R^3;\R^3)$ and $\n(t,\cdot)$ is a normal field to $\Gamma(t)$ for every $t\in[0,T]$. By $\Phi_{-t}\colon \Gamma(t)\to \Gamma_0$ we denote the inverse of $\Phi_t\colon \Gamma_0\to \Gamma(t)$ for $t\in [0,T]$.
 Moreover, we define $\mathbf{D}=\mathbf{D}(t,\z):=D\Phitn(\z)\P_0\colon\R^3\to\R^3$ for every $t\in [0,T], \z\in \Gamma_0$, where $\P_0=\mathbf I-(\n_0\otimes\n_0)$ is the orthogonal projection onto $T_{\z}\gz$, and $\mathbf{D}^{-}=\mathbf{D}^{-}(t,\mathbf x):= D\Phi_{-t}^n(\mathbf x)\P:\R^3\to \R^3$, where $\P=\mathbf I-(\n\otimes\n)$ is the orthogonal projection onto $T_{\mathbf x}\gam(t)$.  Furthermore, we define $\A=\A(t,\z)\colon\R^3\to \R^3$ by
\begin{align}
\A(t,\z):=J^{-1}\mathbf D(t,\z)+\n\otimes \n_0\qquad \text{for all }t\in [0,T], \z\in \Gamma_0,
\end{align}
which is invertible with inverse $\A^{-1}(t,\z):=J\mathbf{D}^{-}+\n_0\otimes \n$.
Then we have
\begin{align}
\norm{J}_{C^4(\mathcal{S}_0)}+\norm{J^{-1}}_{C^4(\mathcal{S})}+\norm{\mathbf D}_{C^4(\mathcal{S}_0)}+\norm{\mathbf{D}^{-}}_{C^4(\mathcal{S})}+\norm{\A}_{C^4(\mathcal{S}_0)}+\norm{\A^{-1}}_{C^4(\mathcal{S})}\leq C(T).
    \label{regularity_basic}
\end{align}
 We now denote by $\nabla_\Gamma$ the tangential gradient, which is for a vector field defined as the covariant derivative, i.e., $\nabla_\Gamma:=\P \nabla \P$, where $\nabla$  the gradient of the function defined on $\Gamma$ and extended in an open neighborhood $U\subseteq \R^3$ of $\Gamma$. Given a differentiable vector field $\v\colon S\to \R^3$, we also define $(\widehat{\nabla}_\gam\v)_{ij}:=(\nablag v_i)_j$, $i,j=1,2,3$. Then $\nablag \v=\P\widehat{\nabla}_\gam\v$. For of a differentiable vector field $\v\colon S\to \R^3$, the tangential divergence is then defined as $\divg \v:=\mathrm{tr}(\nablag \v )$, whereas for a matrix-valued differentiable $\A\colon S\to \R^{3\times 3}$ we define $\divg \A:=\{\divg(\mathbf e_i^T\A)\}_{i=1,2,3}$, where $\{\mathbf e_i\}_i$ is the canonical basis of $\R^3$. 

 We note that $\H:=\nablag\n$ and $H=\text{tr}\H$. Thus
\begin{align}
\norm{H}_{C^4(\mathcal{S})}+\norm{\vn}_{C^4(\mathcal{S})}\leq C(T).
    \label{regularity_basic2}
\end{align}
Let us recall the definitions of  pullback and pushforward maps. For  $\f\colon \mathcal{S}_0\to \R^q$, for some $q\in \N$,
$$
(\tphi_{t}\f)(t,\mathbf x)=\f(t,\Phi_{-t}^n(\mathbf x)),
$$
and  for any $\mathbf g\colon \mathcal{S}\to \R^q$,
$$
(\tphimt\mathbf g)(t,\z):=\mathbf g(t,\Phitn(\z)).
$$
Furthermore surface Piola transform of $\f\colon \mathcal{S}_0\to \R^3$ is defined by
$$
(\phi_{t}\f)(t,\mathbf x)=\A(t,\Phi_{-t}^n(\mathbf x))\f(t,\Phi_{-t}^n(\mathbf x))=\A(t,\Phi_{-t}^n(\mathbf x))(\tphi_t\f)(t,\mathbf x),
$$
as well as its inverse is defined for $\mathbf g\colon \mathcal{S}\to \R^3$ by
$$
(\phi_{-t}\mathbf g)(t,\z):=\A^{-1}(t,\Phitn(\z))\mathbf g(t,\Phitn(\z))=\A^{-1}(t,\Phitn(\z))(\tphimt\g)(t,\z).
$$
We note that the surface Piola transform maps  tangential vector fields to tangential vector field. Moreover, it has the property that $\text{div}_{\gam(t)}\f=0$ for almost any $\mathbf x\in \gam(t)$ if and only if $\text{div}_{\gz}\phi_{-t}\f=0$  for almost any $\z\in \gz$.

We refer to \cite[Section~3]{AGP1} for further details and references. 

\subsection{Evolving Banach and Hilbert spaces}
We use the same notation for function spaces on evolving surfaces as in \cite{AGP1}, which is based on \cite{Def,AE}. For the convenience of the reader we briefly recall important definitions.
For any $t\in[0,T]$ we denote by $W^{k,p}\gt$%
 , where $k\in \mathbb{N}$ and $1\leq p\leq \infty $, the standard Sobolev space and its norm by $\Vert \cdot
 \Vert _{W^{k,p}\gt}$ and as usual $H^{k}\gt=W^{k,2}\gt$. Furthermore, we set $H^{-1}\gt=(H^1\gt)'$.

 For a vector space $X$ of functions defined on $\Gamma(t)$, the space $\mathbf{X}$ is the generic space of tangential vector field or matrix-valued funtions such that each component belongs to $X$. 
 Moreover, we denote by $(\cdot,\cdot )$ the inner product in $\mathbf{L}^{2}\gt$ and by $\Vert \cdot \Vert $
 the induced norm. The inner product of a Hilbert space $H$ is denoted by $(\cdot ,\cdot )_{H}$ and $\Vert \cdot
 \Vert _{H}$ is its induced norm. Similarly as in \cite{V2}, we introduce spaces of solenoidal vector fields
 \begin{align*}
     &\L^2_\sigma\gt:=\{\v\in \L^2\gt: \divg \v=0\text{ a.e. on }\gam(t)\},
     \\
     &\H^1_\sigma\gt:=\{\v\in \H^1\gt: \divg \v=0\text{ a.e. on }\gam(t)\}.
 \end{align*}
 Notice that it holds the Hilbert triplet $\H^1_\sigma\gt\hookrightarrow \L^2_\sigma\gt \hookrightarrow \H^1_\sigma\gt'$.

 For a Banach space $X_0$ (independent of $t$) $L^q (0,T; X_0 )$ is the Bochner space of $X_0$-valued $q$-integrable (or essentially bounded) strongly measurable functions.
 The first order $X_0$-valued one-dimensional Sobolev space is denoted by  $W^{1,p} (0, T ; X_0)$, where $1 \leq p < \infty$ and $H^1 (0, T ; X_0) =W^{1,2}(0,T;X_0)$. Furthermore 
 $C^{0,\alpha} (0,T ; X)$ for $\alpha\in[0,1]$ is the Banach space of all H\"{o}lder continuous functions $f \colon [0,T]\to X_0$ with exponent $\alpha$. 

 Now, if  $X(t)$ is a generic scalar Banach space of functions over $\Gamma(t)$ (and $\mathbf X(t)$ a space of tangential vector fields over $\Gamma(t)$), we have for mappings
\begin{align*}
u\colon [0, T]&\to \bigcup_{t\in[0,T]}\{t\}\times X(t), \quad t\mapsto (t,\bar{u}(t))
\end{align*}
and $1\leq p \leq \infty$ that $u\in L^p_{X}$ if and only if  $t\mapsto \tphi_{-t}\bar{u}(t)\in L^p(0, T; X_0)$ with the norm
 \begin{align*}
     \|u\|_{L^p_X} := \begin{cases}
         \left( \int_0^T \|u(t)\|_{X(t)}^p dt\right)^{1/p} &\text{ if } p<\infty, \\
         \text{ess sup}_{t\in [0,T]} \|u(t)\|_{X(t)} &\text{ if } p=\infty.
     \end{cases}
 \end{align*}
 Usually we identify $u(t)$ and $\bar{u}(t)$.
 If $X(t)=H(t)$ are Hilbert spaces, we equip $L^2_H$ with the inner product
\begin{align*}
(u,v)_{L^2_{H}} &:= \int_0^T (u(t), v(t))_{H(t)}dt. 
\end{align*} 
For the spaces of tangential vector fields $\mathbf X(t)=\L^2_\sigma\gt, \H^1_\sigma\gt,\H^2\gt$, we adopt as pullback the map $\phi_t$ and define, for $p\in [1,\infty]$, $\v\in L^p_{\mathbf X}$ if  $t\mapsto \phi_{-(\cdot)}\bar{\v}(\cdot)\in L^p(0, T; \mathbf X_0)$ with the norm
 \begin{align*}
     \|\v\|_{L^p_\mathbf X} := \begin{cases}
         \left( \int_0^T \|\v(t)\|_{\mathbf X(t)}^p dt\right)^{1/p} &\text{ if } p<\infty, \\
         \text{ess sup}_{t\in [0,T]} \|\v(t)\|_{\mathbf X(t)} &\text{ if } p=\infty.
     \end{cases}
 \end{align*}
 Furthermore, we have $u\in C^1_{X}$ if and only if $t\mapsto \tphi_{-t}u(t)\in C^1([0,T]; X_0)$ and we define its time derivative as 
\begin{align}
\partial^{\bullet}_t u(t)=\dt u(t)= \tphi_t \dfrac{d}{dt} \tphi_{-t} u(t).\label{dtclass}
\end{align}
 These spaces are equipped with the canonical norm
\begin{align*}
\|u\|_{C^1_{X}} := \sup_{t\in [0,T]} \|u(t)\|_{X(t)} +  \sup_{t\in [0,T]} \|\dt u(t)\|_{X(t)}.
\end{align*} 
As in \cite{V2} we define for spaces of divergence-free vector fields $\mathbf X(t)$ as $\L^2_\sigma\gt$ and $\H^1_\sigma\gt$ two time derivatives. The standard time derivative is defined as in \eqref{dtclass}, i.e.,
$\v\in C^1_{\mathbf X}$ if and only if $t\mapsto \tphi_{-t}u(t)\in C^1([0,T]; \X_0)$ and 
\begin{align}
\partial^{\bullet}_t \v(t):=\dt \v(t):= \tphi_t \dfrac{d}{dt} \tphi_{-t} \v(t).\label{dtclass1}
\end{align}
Note that this time derivative is not divergence-free in general. Therefore we also define the time derivative
 \begin{align*}
     \dts \v:=\phi_t\ddt\phi_{-t}\v(t),
 \end{align*}
 which is also consistent with the abstract framework developed in \cite{Def} when applied to the spaces $\L^2_\sigma\gt$ and $\H^1_\sigma\gt$ since it preserves divergence-free tangential vector fields.
 Because of \cite[Lemma 3.6]{V2}, we have the relation
 \begin{align}
\dt\v=\dts\v-\A(\dt\A^{-1})\v,\quad \P\dt\v=\dts\v-\P\A(\dt\A^{-1})\v.
    \label{relation_base1}
\end{align}
These definitions give rise to two notions of weak time derivatives: For $u\in L^2_{H^1}$ a function $v\in L^2_{H^{-1}}$ is the weak time derivative of $u$, (denoted by $v=\dt u$) if 
\begin{align*}
&\int_0^T \langle v(t), \eta(t)\rangle_{H^{-1}(\Gamma(t)), H^1(\Gamma(t))}  = -\int_0^T (u(t),\dt\eta(t)) - \int_0^T \int_{\Gamma(t)} u(t)\eta(t) v_\n H\ds
\end{align*}
for every $\eta\in \mathcal{D}_{H^1}$ (where $\eta\in \mathcal{D}_{X}$ if and only if $t\mapsto \phi_{-t}u(t)\in C_0^\infty((0,T); X_0)$).
As usual the weak derivative $\dt u$ coicides with classical derivative introduced before if $u$ is continously differentiable. 
The same definition holds for the time derivative $\dt\v$ of a   divergence-free tangential vector field $\v$: we have $\w=\dt\u$ in the weak sense if and only if 
\begin{align*}
&\int_0^T \langle \w(t), \boldsymbol\eta(t)\rangle_{\H^{1}_\sigma(\Gamma(t))',\H^1_\sigma(\Gamma(t))} = -\int_0^T (\v(t),\dt\boldsymbol\eta(t)) - \int_0^T \int_{\Gamma(t)} \v(t)\cd\boldsymbol\eta(t) v_\n H\ds
\end{align*}
for any $\boldsymbol\eta\in \mathcal{D}_{{\H^1_\sigma}}$. Finally, because of \eqref{relation_base1} for smooth solutions, we define 
$\tilde{\w}=\dts\u$  if and only if  
\begin{align*}
&\int_0^T \langle \tilde{\w}(t), \boldsymbol\eta(t)\rangle_{\H^{1}_\sigma(\Gamma(t))',\H^1_\sigma\gt}\\& = -\int_0^T (\v(t),\dts\boldsymbol\eta(t)) - \int_0^T \int_{\Gamma(t)} \v(t)\boldsymbol\eta(t) v_\n H\ds+\int_0^T\ig (\A\dt(\A^{-1})+(\A\dt(\A^{-1}))^T)\v\cd \boldsymbol\eta\ds
\end{align*}
for every $\boldsymbol\eta\in \mathcal{D}_{{\H^1_\sigma}}$. Then \eqref{relation_base1} also holds for the weak notions of time derivatives.

We then define the Banach spaces
\begin{align*}
H^1_{{\H^{1}_\sigma}'} := \left\{ \v\in L^2_{\H^{1}_\sigma} \colon \dts \v\in L^2_{{\H^{1}_\sigma}'}\right\} \, \text{ with } \,\, \|\v\|_{H^1_{{\H^{1}_\sigma}'}} := \|\v\|_{L^2_{{\H^{1}_\sigma}}} + \| \dts \v \|_{L^2_{{\H^{1}_\sigma}'}}+\|\v(0)\|,
\end{align*}
and 
\begin{align*}
H^1_{H^{-1}} := \left\{ u\in L^2_{H^1} \colon \dt u\in L^2_{H^{-1}}\right\} \, \text{ with } \,\, \|u\|_{H^1_{H^{-1}}} := \|u\|_{L^2_{H^1}} + \| \dt u \|_{L^2_{H^{-1}}}+\|u(0)\|.
\end{align*}
 In the space $H^1_{{\H^{1}_\sigma}'}$ one could also adopt the weak notion of the derivative $\dt$, but in this case, since $\dt\v$ is not divergence-free, one should ask for $\dt\v_{|\H^1_{\sigma}}\in L^2_{{\H^1_\sigma}'}$. Notice that knowing only the action of $\dt\v$ on the functions in $\H^1_\sigma(\Gamma(t))$ does not allow in principle to reconstruct the entire time derivative $\dt\v$, which has to be obtained from \eqref{relation_base1}. On the other hand, if we know that $\dts\v\in L^2_{{\H^1_\sigma}'}$, then from \eqref{relation_base1} is immediate to infer $\dt\v\in L^2_{{\H^1}'}$.  We also define $H^1_{H^1}$, which denotes the space of those $u\in L^2_{H^1}$ which have a more regular weak time derivative $\dt u\in L^2_{H^1}$. In general we denote by $W_{X,Y}$ the spaces of functions $u$ such that $u\in L^2_X$ and $\dt u\in L^2_Y$ when $X$ and $Y$ do not coincide, whereas we write $H^1_{X}:=W_{X,X}$.
Furthermore, $W(\H^1_\sigma\cap \H^2,{\L^2_\sigma})$ denotes the space of those $\v\in L^2_{\H^1_\sigma\cap \H^2}$ which have a more regular weak time derivative $\dts \v\in L^2_{\L^2_\sigma}$. Note that we prefer to define the space in this way, since if we consider in the definition of the last space the derivative $\dt\v$, this quantity would not be divergence-free nor tangential, whereas, as it is standard in the fixed domain case, we expect a divergence-free velocity time derivative to belong to $\L^2_\sigma\gt$. Observe that it can be proven, exploiting the regularity of \eqref{regularity_basic}, that there is an evolving space equivalence in the sense of \cite[Definition 2.31]{Def} between the spaces $H^1_{H^{-1}}$ and $H^1(0,T;H^{-1}(\gz))\cap L^2(0,T;H^1(\gz))$, as well as between $W(H^3,H^{-1})$ and $H^1(0,T;H^{-1}(\gz)\cap L^2(0,T;H^3(\gz))$, between $W(H^3,H^{1})$ and $H^1(0,T;H^{1}(\gz)\cap L^2(0,T;H^3(\gz))$, and between $W({H^{1},H^1})$ and $H^1(0,T;H^{1}(\gz))$. Analogously, by arguing for instance, in a similar way as in \cite[Lemma 3.7]{V2}, there is an evolving space equivalence between the spaces $W(\H^2\cap \H^1_\sigma, \L^2_\sigma)$ and $H^1(0,T;\L^2_\sigma(\gz))\cap L^2(0,T;\H^2(\gz))$. Clearly, the same evolving space equivalence holds between $W(\H^2, \L^2)$ and $H^1(0,T;\L^2(\gz))\cap L^2(0,T;\H^2(\gz))$.
Therefore, by standard embeddings for the Bochner spaces over $\gz$ and recalling these evolving space equivalences, we deduce the following 
\begin{lemma}
\label{continuous_embeddings}
    The following embeddings hold 
    \begin{align*}
    \non& H^1_{H^{-1}}\hookrightarrow C^0_{L^2},\quad H^1_{H^1}\hookrightarrow C^0_{H^1},\quad W_{H^3,H^{-1}}\hookrightarrow C^0_{H^1},\\&
    W_{H^3,H^1}\hookrightarrow C^0_{H^2},\quad W(\H^2\cap \H^1_\sigma, \L^2_\sigma)\hookrightarrow C^0_{\H^1_\sigma},\quad W(\H^2, \L^2)\hookrightarrow C^0_{\H^1}.
    \end{align*}
    
\end{lemma}

Finally, if  functions are not defined on the entire interval of time $[0,T]$, but only on $[0,T_1]$, with $0<T_1\leq T$, the corresponding spaces are denoted by an extra $(T_1)$ as for instance, $L^p_{X(T_1)}$. We note that Lemma \ref{continuous_embeddings} also holds for the spaces on $[0,T_1]$ and the embedding constants only depend on $T$ if the norms of the initial value is included in the norm appropriately. 
\subsection{Some Uniform Estimates on $[0,T]$}
In this section we present some estimates which also hold uniformly in time on $[0,T]$. In particular, we have the following Poincaré inequality (see also \cite{Aubin}): there exists a constant $C_P(T)>0$ such that, for any $u\in H^1(\Gamma(t))$,
$$
\norm{u(t)-\overline{u}(t)}\leq C_P\norm{\nabla_{\Gamma}u},\quad \forall t\in[0,T],
$$
where $\overline{u}(t):=\dfrac{\intg u(t)\ds}{\norma{\gam(t)}}$.  
Furthermore, we have the following adaptation of Korn's inequality (see \cite[Lemma 3.3]{V2}): there exists $C=C(T)>0$ such that, for any $\v\in \H^1\gt$, it holds
\begin{align}
\norm{\v}_{\H^1\gt}\leq C(T)(\norm{\v}+\norm{\E_S(\v)}),\quad \text{ for a.a. }t\in[0,T].
\label{Korn}
\end{align}
Then, have the following Sobolev-Gagliardo-Nirenberg's inequality:
\begin{align}
\norm{f}_{L^p\gt}\leq C(T)\sqrt p\norm{f}^\frac2p\norm{f}_{H^1\gt}^{1-\frac2p},\quad \forall p\in[2,+\infty),\quad\forall f\in H^1\gt,\quad \forall t\in[0,T].
\label{Gagliardo}
\end{align}
\begin{proof}
It is well known  that, for any sufficiently regular open set  $\Omega\subset \R^2$ it holds
$$
\norm{f}_{L^p(\Omega)}\leq C(\Omega)\sqrt p\norm{f}^\frac2p_{L^2(\Omega)}\norm{f}_{H^1(\Omega)}^{1-\frac2p},\quad \forall p\in[2,+\infty),\quad \forall f\in H^1(\Omega).
$$
 It then holds, by a standard partition of unity argument applied to $\Gamma_0$, that there exists $C(\gz)>0$ such that
\begin{align*}
\norm{f}_{L^p(\gz)}\leq C(\gz)\sqrt p\norm{f}^\frac2p_{L^2(\gz)}\norm{f}_{H^1(\gz)}^{1-\frac2p},\quad \forall p\in[2,+\infty),\quad \forall f\in H^1(\Gamma_0).
\end{align*}
Thanks to the regularity of the flow map \eqref{regularity_basic} and the compatibility of the space involved, we know that 
$$
\norm{f}_{L^p\gt}\leq C(T)\norm{\tphimt f}_{L^p(\gz)}\leq C(T)\sqrt p\norm{\tphimt f}^\frac2p_{L^2(\gz)}\norm{\tphimt f}_{H^1(\gz)}^{1-\frac2p}\leq C(T)\sqrt p\norm{f}^\frac2p\norm{ f}_{H^1\gt}^{1-\frac2p},
$$
concluding the proof.
\end{proof}
Clearly (see also \cite[Lemma 3.4]{V2}), the inequality \eqref{Gagliardo} also follows for any (tangential) vector $\f\in \H^1\gt$ by considering it component-wise.  By similar arguments we also deduce the validity of a uniform Agmon's inequality (clearly also valid for tangential vector fields):
\begin{align}
\norm{f}_{L^\infty\gt}\leq C(T)\norm{f}^\frac12\norm{f}_{H^2\gt}^\frac12\quad \forall f\in H^2\gt,\quad \forall t\in[0,T].
    \label{Agmon}
\end{align}
In conclusion, let us assume that $f(t)\in H^2\gt$, for almost any $t\in[0,T]$, satisfies, for some $\tilde{\omega}>0$, 
$$
-\Delta_\gam f(t)+\tilde{\omega}f(t)=h(t).
$$
Then we can show that it holds
\begin{align}
\norm{f(t)}_{H^{k+2}\gt}\leq C(T)\norm{h(t)}_{H^k\gt},\quad k=1,2,\quad\text{ for a.a. } t\in[0,T],
\label{uniformelliptic}
\end{align}
where $C(T)>0$ does not depend on $t\in[0,T]$. The proof of this result exploits the elliptic regularity theory applied on each surface $\gam(t)$ as well as the regularity \eqref{regularity_basic} on the flow map. For the sake of brevity we only refer to the proof of \cite[Lemma 7.4]{AGP1} for a very similar argument.

\subsection{Results on Time Differentiation of some Bilinear Forms}


We have the following essential theorem to deal with some variations on the classical Transport Theorem: 
\begin{theorem}
Let the assumptions on the flow map $\Phitn$ of Section \ref{regflowmap1} hold. 
\begin{itemize}

\item For any $\eta, \phi\in H^1_{L^2}$ it holds
\begin{align}
&\frac{d}{dt}\intg \eta\phi\ds=\intg \dt \eta \phi\ds +\intg \eta \dt\phi\ds+\intg  \eta\phi\ \divg\V_\n\ds. 
\label{dt1}
\end{align}
\item For any $\eta,\phi\in H^1_{H^1}$ it holds 
\begin{align}
\label{dt3}&\frac{d}{dt}\intg \nablag\eta\cdot \nablag\phi\ds \\&=\non\intg \nablag\dt \eta \cdot\nablag \phi\ds +\intg \nablag\eta\cdot \nablag\dt\phi\ds+\intg  \nablag\eta\cdot \nablag\phi\ \divg\V_\n\ds-2\intg \nablag\phi\cdot\E_S(\V_\n)\nablag\eta\ds.
\end{align}
\item For any $f\in H^1_{H^1}$ and any $\g\in H^1_{\L^2}$ ($\g$ tangential vector field), it holds
\begin{align}
\non&\frac{d}{dt}\intg \nablag f\cdot \g\ds\\&=\intg  \nablag f\cdot \dt \g\ds+\intg \nablag\dt f\cdot \g\ds-\intg \g\cdot (\nablag^T\V_n)\nablag f\ds+\intg \nablag f\cdot \g\ \divg\V_\n\ds.
\label{dt2}
\end{align}

\item For any $\u,\v\in H^1_{\H^1}$ (tangential vector fields) and any function $\eta\in L^\infty_{L^\infty}$ such that $\dt\eta\in L^\infty_{L^\infty}$, it holds
\begin{align}
\nonumber\frac{d}{dt}\intg \eta \E_S(\u):\E_S(\v)\ds&=\int_{\gam}\dt\eta\E_S(\u):\E_S(\v)\ds\\&\nonumber+\int_\gam \eta\E_S(\dt \u):\E_S(\v)\ds+\int_\gam \eta\E_S( \u):\E_S(\dt\v)\ds\\&\nonumber+2\int_{\gam}\eta\mathcal{S}((\n\otimes \n)\widehat{\nabla}_\gam\V_\n)\widehat{\E}_S(\u):\E_S(\v)\ds+2\int_\gam\eta\widehat{\E}_S(\u)\mathcal{S}((\n\otimes \n)\widehat{\nabla}_\gam\V_\n):\E_S(\v)\ds\\&\nonumber+2\int_{\gam}\eta{\E}_S(\u):\mathcal{S}((\n\otimes \n)\widehat{\nabla}_\gam\V_\n)\widehat{\E}_S(\v)\ds+2\int_{\gam}\eta{\E}_S(\u):\widehat{\E}_S(\v)\mathcal{S}((\n\otimes \n)\widehat{\nabla}_\gam\V_\n)\ds\\&-\int_\gam \eta\mathcal{S}((\nablag\u)\nablag\V_\n):\E_S(\v)\ds-\int_\gam \eta\E_S(\u):\mathcal{S}((\nablag\v)\nablag\V_\n)\ds\nonumber\\&+\intg \eta \E_S(\u):\E_S(\v)\ \divg\V_\n\ds,
\label{Ess}
\end{align}
where $\widehat{\E}_{S} (\cdot):=\frac{\widehat{\nabla}_{\gam}\cdot+\widehat{\nabla}_\Gamma^T\cdot}{2}$ and $\mathcal{S}(\A)$ is the symmetric part of a matrix $\A$. Here by the matrix $(\nablag\v)\nablag\V_\n$ we denote the matrix product between $\nablag\v$ and $\nablag\V_\n$.

\end{itemize}
\end{theorem}
\begin{remark}
We observe that, by a standard density argument, \eqref{Ess} also holds if $\u,\v\in L^2_{\H^2}\cap H^1_{\L^2}$ and $\eta\in L^\infty_{W^{1,\infty}}$ with $\dt\eta\in L^2_{H^1}$, up to substituting the second line of the right-hand side by 
$$
-\intg \divg(\eta\E_S(\u))\cdot \dt\v\ds-\intg \divg(\eta\E_S(\v))\cdot \dt\u\ds,
$$
as it can be easily seen integrating by parts the corresponding terms.
\label{basic_rem}
\end{remark}
\begin{proof}
The proofs of \eqref{dt1}-\eqref{dt3} can be found in \cite[Section 8.2]{DzEll}, whereas the proof of \eqref{dt2} is shown, for instance, in \cite[Appendix A.1]{CEGP}. We are left to prove \eqref{Ess}. To this aim, we first observe that, under the same assumptions of \eqref{dt2}, it holds, applying \eqref{dt1} component-wise, 
\begin{align*}
&\frac{d}{dt}\intg \nablag f\cdot \g\ds
 \\&=\intg \dt \nablag f\cdot \g \ds+\intg  \nablag f\cdot \dt\g\ds+\intg   \nablag f\cdot \g\ \divg\V_\n\ds,
\end{align*}
so that, by comparison with \eqref{dt2} and by the arbitrariness of the tangential vectors $\g\in H^1_{\L^2}$, we obtain
\begin{align*}
\P\dt\nablag f=\nablag\dt f-\nablag^T\V_n\nablag f.
\end{align*}
Observe that this identity can be deduced also from \cite[Lemma 2.6]{Dzk}.
If we now consider a tangential vector field $\u\in H^1_{\H^1}$, we then have (by choosing for instance $f=\u_j$, $j=1,2,3$, in the identity above and after some standard computation) 
\begin{align*}
\P(\dt\widehat{\nabla}_{\gam} \u)\P=\nablag\dt \u-(\nablag\u)\nablag\V_\n,
\end{align*}
so that, by taking the transpose and summing up dividing by two,
\begin{align}
\P(\dt\widehat{\E}_{S} (\u))\P=\E_S(\dt \u)-\mathcal{S}((\nablag\u)\nablag\V_\n).
\label{Es1}
\end{align}
We now recall that, as in \cite[Lemma 2.2]{V1}, it holds 
$$
\dt\P=2\mathcal{S}((\n\otimes \n)\widehat{\nabla}_\gam\V_\n).
$$
By applying \eqref{dt1} (if we think component-wise), we get, for $\u,\v\in H^1_{\H^1}$ and $\eta\in L^\infty_{L^\infty}$ such that $\dt\eta\in L^\infty_{L^\infty}$, 
\begin{align}
\nonumber&\frac d{dt}\intg \eta \E_S(\u):\E_S(\v)\ds\\&=\int_{\gam}\dt\eta\E_S(\u):\E_S(\v)\ds+ \int_{\gam}\eta\dt\E_S(\u):\E_S(\v)     & \nonumber\\&+\int_{\gam}\eta\E_S(\u):\dt\E_S(\v)\ds+\intg \eta \E_S(\u):\E_S(\v)\ \divg\V_\n\ds.
\label{dtE}
\end{align}
Now, by the product rule for differentiation, we get 
\begin{align*}
\dt\E_S(\u)=\dt (\P\widehat{\E}_S(\u)\P)&=(\dt \P)\widehat{\E}_S(\u)\P+\P\dt\widehat{\E}_S(\u)\P+\P\widehat{\E}_S(\u)\dt\P\\&
=2\mathcal{S}((\n\otimes \n)\widehat{\nabla}_\gam\V_\n)\widehat{\E}_S(\u)\P+\P\dt\widehat{\E}_S(\u)\P+2\P\widehat{\E}_S(\u)\mathcal{S}((\n\otimes \n)\widehat{\nabla}_\gam\V_\n).
\end{align*}
For the middle term in the identity above, observe that it holds
\begin{align*}
\int_{\gam}\eta\P\dt\widehat{\E}_S(\u)\P:\E_S(\v)\ds=\int_{\gam}\eta(\E_S(\dt \u)-\mathcal{S}((\nablag\u)\nablag\V_\n)):\E_S(\v)\ds,
\end{align*}
so that,  since $\P\E_S(\v)\P=\E_S(\v)$,
\begin{align*}
\int_{\gam}\eta\dt\E_S(\u):\E_S(\v)\ds&=2\int_{\gam}\eta\mathcal{S}((\n\otimes \n)\widehat{\nabla}_\gam\V_\n)\widehat{\E}_S(\u):\E_S(\v)\ds\\&+2\int_\gam\eta\widehat{\E}_S(\u)\mathcal{S}((\n\otimes \n)\widehat{\nabla}_\gam\V_\n):\E_S(\v)\ds\\&+\int_\gam \eta\E_S(\dt \u):\E_S(\v)\ds-\int_\gam \eta\mathcal{S}((\nablag\u)\nablag\V_\n):\E_S(\v)\ds
\end{align*}
Therefore, \eqref{Ess} follows from \eqref{dtE} inserting the identity above and the one for  $\int_{\gam}\eta\E_S(\u):\dt\E_S(\v)\ds$, which is identical up to exchanging the roles of $\u$ and $\v$.
\end{proof}
\subsection{Assumptions on the Potential, the Viscosity and the Density}
In the main result of this paper we consider the following assumptions:
\begin{enumerate}[label=(\subscript{H}{{\arabic*}})] 
\item \label{h1}$\Psi\in C^0([-1,1])\cap C^{4}((-1,1))$. Furthermore, the potential $\Psi:[-1,1]\to \R$ can be written as
    \begin{equation*}
        \Psi(s)=F(s)-\frac{\theta_0}{2}s^2 \quad\text{for all $s\in [-1,1]$}
    \end{equation*}
    with a given constant $\theta_0>0$. 
   Moreover $F\in C^0([-1,1])\cap C^{4}((-1,1))$ has the properties 
    \begin{equation*}
    \lim_{r\rightarrow -1}F^{\prime }(r)=-\infty ,
    \quad \lim_{r\rightarrow 1}F^{\prime }(r)=+\infty ,
    \quad F^{\prime \prime }(s)\geq {\theta},
    \end{equation*}
    for all $s\in (-1,1)$ and  for some $0<\theta<\theta_0$.
    Without loss of generality, we further assume $F(0)=0$ and $F'(0)=0$. 
    In particular, this implies that $F(r)\geq 0$ for all $r\in [-1,1]$.

    \noindent We also assume that there exists $\beta>\frac12$ such that
    \begin{equation}
    \frac{1}{F^{\prime }(1-2\delta )}=O\left( \frac{1}{|\ln (\delta )|^{\beta }}%
    \right) ,\quad\text{ }\dfrac{1}{|F^{\prime }(-1+2\delta )|}=O\left( \frac{1}{%
    |\ln (\delta )|^{\beta }}\right)   \label{assmp}
    \end{equation}
    as $\delta\to 0^+$.
    
\item \label{h2}$\nu\in C^{2}(\R)$, such that $\nu(s) \geq \nu_* > 0$  for every $s \in \R$
and some $\nu_* > 0$.
\item \label{rho}
$  \rho(s):=\frac{\widetilde{\rho}_1+\widetilde{\rho}_2}{2}+\frac{\widetilde{\rho}_2-\widetilde{\rho}_1}{2}s,\ \forall s\in \R$, where we set $\rho_*=\min_{s\in[-1,1]}\rho(s)>0$ and $\rho^*=\max_{s\in[-1,1]}\rho(s)>0$.
\end{enumerate}

\begin{remark}
A common form for the
viscosity is the following
\begin{equation*}
\nu (s)=\nu _{1}\frac{1+s}{2}+\nu _{2}\frac{1-s}{2},\quad \forall s\in \lbrack -1,1],
\end{equation*}%
which can be easily extended on the whole $\mathbb{R}$ in such way to comply
with \ref{h2}.
\end{remark}
\begin{remark}\label{REM:LOG}
Notice that the logarithmic potential, which is given by
\begin{equation}
    F(s)=\frac{\theta}{2}((1+s)\text{ln}(1+s)+(1-s)\text{ln}(1-s))
    \quad\text{for all $s\in(-1,1)$},
\label{F:LOG}
\end{equation}
complies with assumption \ref{h1}, with some $0<\theta<\theta_0$. Nevertheless, much more general singular potentials are allowed by \eqref{assmp} (see, for instance, \cite{GalP}).
\end{remark}
\section{Main Result}\label{sec:main}
In this section we present the main result of this paper, i.e., the existence and uniqueness of a global strong solution to \eqref{mainp}. We first give the definition of a strong solution on an interval $[0,T_0]$, with $0<T_0\leq T$, to problem \eqref{mainp}, namely 
\begin{definition}[Strong solution]
	\label{strong}
	Let us set $\mathbf{v}_0 \in \H^1(\gz)$
	and $\vphi_0\in {H}^2 (\gz)$ such that $-\Delta\vphi_0+\Psi'(\vphi_0)=:\mu_0\in H^1(\gz)$, $\norm{\vphi_0}_{L^\infty(\gz)}\leq1$ and  $\norma{\overline{\vphi}_0}=\norma{\frac{\int_{\gz}\vphi_0\ds}{\norma{\gz}}}<1$. We call
$(\mathbf{v},p,\vphi)$ a strong solution to \eqref{mainp} on $[0,T_0]$ if 
	\begin{align*}
	&\mathbf{v} \in   L ^2_{\H^2(T_0)}\cap H^1_{\L^2(T_0)} ,\quad p\in L^2_{H^1(T_0)},\quad 
	\vphi\in L^\infty_{W^{2,q}(T_0)} \cap H^1_{{H}^1(T_0)} ,\quad \forall q\in[2,+\infty),\\&
\norma{\varphi}<1\quad\text{a.e. in }\gam(t) \quad\text{for a.a. }t\in[0,T_0],\quad 
	\mu\in {L}^\infty_{H^1(T_0)}\cap L^2_{H^3(T_0)},
	\end{align*}
	and it satisfies \eqref{mainp} in the almost everywhere sense.

\end{definition} 
\begin{remark}
Since, by Lemma \ref{continuous_embeddings}, $L ^2_{\H^2}\cap H^1_{\L^2}\hookrightarrow C_{\H^1}$, and $H^1_{H^1}\hookrightarrow C_{H^1}$, the initial conditions $\v(0)=\v_0$ and $\vphi(0)=\vphi_0$ make sense.
\end{remark}
Our main result then reads 
\begin{theorem}[Global well-posedness of strong solutions]
	\label{strong1}
	Assume the regularity assumptions on $\Phitn$ stated in Section \ref{regflowmap1}, and assumptions \ref{h1}-\ref{rho} for $\Psi$, $\nu$, and $\rho$. Let $\mathbf{v}_0 \in\H^1(\gz)$
	and $\vphi_0\in {H}^2 (\gz)$, be such that $-\Delta_{\Gamma_0}\vphi_0+\Psi'(\vphi_0)=:\mu_0\in H^1(\gz)$, $\norm{\vphi_0}_{L^\infty(\gz)}\leq1$ and $\norma{\overline{\vphi}_0}<1$. Then there exists a unique global strong solution $(\mathbf{v},p, \vphi)$ to \eqref{mainp} on $[0,T]$ in the sense of Definition \ref{strong}. Furthermore, the following additional regularity holds: there exists $\delta>0$ depending only on the parameters of the problem, $F$, $T$, the initial energy $E_{tot}(\rho_0,\v_0,\vphi_0)$, and $\norm{\mu_0}_{H^1(\gz)}$, such that
 \begin{align}
    \sup_{t\in[0,T]} \norm{\vphi}_{C\gt}\leq 1-\delta. \label{separaz}
 \end{align}
 This entails
 \begin{align}
 \vphi\in L^\infty_{H^3}\cap L^2_{H^4},\quad \mu\in H^1_{H^{-1}}.
     \label{H3}
 \end{align}
\end{theorem}
\begin{remark}
    We observe that, exploiting the estimates obtained in the first part of the proof of Theorem \ref{strong1} (see, in particular, the energy estimate \eqref{etot} below), one can also prove the existence of global-in-time weak solutions to the same problem \eqref{mainp}, under weaker assumptions on the regularity of the intial data. The same result has been obtained, in the case of matched densities (i.e., $\rho=const$) in \cite{ES}. Using our approach, the existence result of weak solutions can also be extended to the more general case of unmatched densities (we also refer to \cite{ADG, AGP} for a  definition of weak solutions in the unmatched densities case in bounded domains). We leave the details to the interested reader.
\end{remark}
Before proving this result, in the next section we analyze the well-posedness of a convective Cahn-Hilliard equation on evolving surfaces. This is a milestone to prove our main result.

\section{Regularity Results for the Advective Cahn-Hilliard Equation on Evolving Surfaces}\label{sec:ConvCH}
In this section we address the well-posedness of the following problem:
\begin{align}
   \begin{cases}
    \dt\vphi-\Delta_{\Gamma}\mu+\nablag\vphi\cdot \v=0,&\quad\text{ on }\Gamma,\\
         \mu=-\Delta_\Gamma\vphi+\Psi'(\vphi),&\quad \text{ on }\Gamma,\\
         \vphi(0)=\vphi_0,&\quad\text{ on }\gz,\end{cases}
         \label{pr}
\end{align}
where $\v=\u+\VVn$ is divergence-free, and $\u$ is a tangent vector.
In particular, it holds the following
\begin{theorem}
	\label{existence}
	Let $\Psi$ satisfy assumption \ref{h1}  and the flow map satisfy the assumptions of Section \ref{regflowmap1}, and let $\vphi_0\in H^2(\gz)$ such that $\norma{\overline{\vphi}_0}<1$, $\norm{\vphi_0}_{L^\infty(\gz)}\leq1$, and $\mu_0:=-\Delta_{\gz}\vphi_0+\Psi'(\vphi_0)\in H^1(\gz)$. Assume also that $\v=\u+\VVn$ and $\u\in L^2_{\H^1(\tT)}$, for some $\tT\in(0,T]$, and $\divg\v\equiv 0$. Then there exists a unique pair $(\vphi,\mu)$ with
\begin{align*}	
 &\vphi\in L^\infty_{H^3(\tT)}\cap L^2_{H^4(\tT)}\cap H^1_{H^{1}(\tT)},\\&
 \mu\in L^\infty_{H^1(\tT)}\cap L^2_{H^3(\tT)}\cap H^1_{H^{-1}(\tT)},\\&
 F'(\vphi)\in L^\infty_{L^{p}(\tT)},\quad \forall p\geq 2
 \end{align*}
 such that it holds
 $$
\sup_{t\in[0,\tT]}\norm{\vphi(t)}_{C\gt}\leq 1-\delta,
 $$
 for some $\delta=\delta(T,\norm{\u}_{L^2_{\H^1(\tT)
 }})\in(0,1)$, continuously decreasing towards 0 with the increase of the $L^2_{\H^1(\tT)}$ norm of $\u$. The pair $(\vphi,\mu)$ is a strong solution to problem \eqref{pr}, i.e., it satisfies the equation in the almost everywhere sense (for almost any $t\in[0,T]$ almost everywhere on $\Gamma(t)$). In particular, there exists $C(T,\norm{\u}_{L^2_{\H^1(\tT)
 }})>0$, continuously increasing with the $L^2_{\H^1(\tT)}$ norm of $\u$, such that
 \begin{align} \norm{\mu}_{L^\infty_{H^1(\tT)}}+\norm{\mu}_{L^2_{H^3(\tT)}}+\norm{\mu}_{H^1_{H^{-1}(\tT)}}+\norm{\vphi}_{L^\infty_{H^3(\tT)}}+\norm{\vphi}_{H^1_{H^1(\tT)}}\leq C(T,\norm{\u}_{L^2_{\H^1(\tT)
 }}).
 \end{align}

\end{theorem}
\begin{remark}
    Note that, since $\vphi\in H^1_{H^1(\tT)}\hookrightarrow C_{H^1(\tT)}$, the initial condition $\vphi(0)=\vphi_0$ makes sense.
\end{remark}
\begin{remark}
  Due to the assumption $\mu_0\in H^1(\Gamma_0)$, there exists $\delta_0>0$ such that 
    $$
\norm{\vphi_0}_{C^0(\gz)}\leq 1-\delta_0.
    $$
    Indeed, it is immediate to see that we have, thanks to \eqref{Gagliardo},  
    $$
\norm{\mu_0}_{L^p(\gz)}+\norm{F'(\vphi_0)}_{L^p(\gz)}\leq C\sqrt{p},    \quad \forall p\geq 2.
    $$
    Then we can repeat the proof of \cite[Theorem 3.3]{GalP} for the time $t=0$ and deduce the result (see also the proof of \eqref{sepp} below for a similar argument). 
    Therefore, since $\mu_0\in H^1(\gz)$ and, exploiting the separation property above, $\Psi'(\vphi_0)\in H^1(\gz)$, it is immediate to deduce by elliptic regularity that $\vphi_0\in H^3(\gz)$.
    \label{necessaryrem}
\end{remark}
\begin{proof}
\textbf{Uniqueness.} First of all, uniqueness also holds for weak solutions and the proof is an immediate adaptation, for instance, of \cite{DE}. Let us introduce the inverse Laplacian operator $\mathcal{N}(t)$ such that, for any $f\in H^{-1}\gt$ with $<f,1>_{H^{-1},H^1}=0$, it holds (recall that $\mathcal{N}(t)f\in H^1(\Gamma(t))$)
\begin{align*}
    \intg \nablag \mathcal{N}(t)f\cdot \nablag\eta\ds=<f,\eta>_{H^{-1},H^1},\quad \forall \eta \in H^1\gt.
\end{align*}
This operator defines a norm on $H^{-1}\gt$ with zero mean (see for instance, \cite[Appendix C]{DE}) which we denote as $\norm{f}_{\sharp}=\norm{\nabla\mathcal{N}(t)f}$. Clearly this norm changes with $t\in[0,T]$. Furthermore, it holds, by standard elliptic regularity (see, e.g., \eqref{uniformelliptic}), 
\begin{align}
\norm{\mathcal{N}(t)f}_{H^2\gt}\leq C(T)\norm{f}, \quad \forall f\in L^2\gt\quad\text{s.t. }\overline{f}=0.
    \label{ctt}
\end{align}
We then consider two initial data $\vphi_{0,1}$ and $\vphi_{0,2}$ satisfying the assumptions of the theorem and such that $\overline{\vphi}_{0,1}=\overline{\vphi}_{0,2}$. Let  us define as $(\vphi_i,\mu_i)$, $i=1,2$, the corresponding solutions. By conservation of total mass, this means that $\overline{\vphi}_1\equiv \overline{\vphi}_2$. Let us then consider the variables $\vphi:=\vphi_1-\vphi_2$ and $\mu:=\mu_1-\mu_2$. They satisfy
 \begin{align}
     \dt \vphi+\v\cdot \nablag\vphi-\Delta_\gam\mu=0.
     \label{uni1}
 \end{align}
 Now we recall that, by \eqref{dt1},

 \begin{align*}
    \frac12 \frac{d}{dt}\norm{\nablag\mathcal{N}(t)\vphi}^2&=\intg \nablag\mathcal{N}(t)\vphi\cdot \nablag\dt \mathcal{N}(t)\vphi\ds\\&\non+\frac12\intg\norma{\nablag\mathcal{N}(t)\vphi}^2 Hv_\n\ds-\frac12\intg \nablag\mathcal{N}(t)\vphi\cdot \E_S(\VVn)\nablag\mathcal{N}(t)\vphi\ds\\&=\intg \vphi\dt \mathcal{N}(t)\vphi\ds+\frac12\intg\norma{\nablag\mathcal{N}(t)\vphi}^2 Hv_\n\ds-\frac12\intg \nablag\mathcal{N}(t)\vphi\cdot \E_S(\VVn)\nablag\mathcal{N}(t)\vphi\ds\\&
    =\frac d{dt}\intg \vphi\mathcal{N}(t)\vphi\ds-\intg \dt\vphi\mathcal N(t)\vphi\ds-\intg \vphi\mathcal{N}(t)\vphi Hv_\n\ds\\&+\frac12\intg\norma{\nablag\mathcal{N}(t)\vphi}^2 Hv_\n\ds-\frac12\intg \nablag\mathcal{N}(t)\vphi\cdot \E_S(\VVn)\nablag\mathcal{N}(t)\vphi\ds,
 \end{align*}
 so that, since $\norm{\vphi}_\sharp^2=\intg \vphi\mathcal{N}(t)\vphi\ds=\norm{\nablag\mathcal{N}(t)\vphi}^2$, we infer
\begin{align*}
   &\frac12 \frac d{dt}\norm{\vphi}_\sharp^2=\intg \dt\vphi\mathcal N(t)\vphi\ds+\intg \vphi\mathcal{N}(t)\vphi Hv_\n\ds\\&-\frac12\intg\norma{\nablag\mathcal{N}(t)\vphi}^2 Hv_\n\ds+\frac12\intg \nablag\mathcal{N}(t)\vphi\cdot \E_S(\VVn)\nablag\mathcal{N}(t)\vphi\ds
\end{align*}
 By testing equation \eqref{uni1} with $\eta=\mathcal{N}(t)\vphi$, we then deduce
\begin{align}
     &\nonumber\frac12 \frac d{dt}\norm{\vphi}_\sharp^2+(\nablag\mu,\nablag\mathcal{N}(t)\vphi)+(\v\cdot \nablag\vphi,\mathcal{N}(t)\vphi)\\&=\intg \vphi\mathcal{N}(t)\vphi Hv_\n\ds-\frac12\intg\norma{\nablag\mathcal{N}(t)\vphi}^2 Hv_\n\ds+\frac12\intg \nablag\mathcal{N}(t)\vphi\cdot \E_S(\VVn)\nablag\mathcal{N}(t)\vphi\ds.
     \label{Gron}
\end{align}
 First, we observe that
 \begin{align*}
     (\nablag\mu,\nablag\mathcal{N}(t)\vphi)=(\mu,\vphi)=\intg (F'(\vphi_1)-F'(\vphi_2))\vphi-\theta_0\norm{\vphi}^2+\norm{\nablag\vphi}^2\ds,
 \end{align*}
 and note that 
 $$
 \intg (F'(\vphi_1)-F'(\vphi_2))\vphi\ds\geq \theta \norm{\vphi}^2.
 $$
 Furthermore, by Cauchy-Schwarz and Young's inequalities,
 \begin{align}
     \theta_0\norm{\vphi}^2=\theta_0\intg \nablag\mathcal{N}(t)\vphi\cdot \nablag\vphi\ds
\leq C\norm{\vphi}_\sharp^2+\frac{1}{4}\norm{\nablag\vphi}^2.
\label{interpola}
 \end{align}
 Then, concerning the advective term, by H\"{o}lder's, Young's and Sobolev-Gagliardo-Nirenberg's inequalities (see also \eqref{Gagliardo}), after integrating by parts, applying Poincaré's inequality (recall that $\overline{\mathcal N\vphi}\equiv0$) and recalling $\v=\VVn+\u$, \eqref{ctt}, and the same interpolation as in \eqref{interpola},
 \begin{align*}
     \norma{(\v\cdot \nablag\vphi,\mathcal{N}(t)\vphi)}&\leq\norma{(\v\vphi,\nablag\mathcal{N}(t)\vphi)}+\norma{\intg \vphi\mathcal{N}(t)\vphi Hv_\n\ds}\\&\leq 
C(1+\norm{\u})\norm{\vphi}_{L^4\gt}\norm{\nablag\mathcal{N}(t)\vphi}_{\L^4\gt}+C\norm{\vphi}\norm{\nablag\NN\vphi}\\&
\leq C(1+\norm{\u})\norm{\vphi}^\frac12\norm{\nablag\vphi}^\frac12\norm{\vphi}_{\sharp}^\frac12\norm{\vphi}^\frac12+C\norm{\vphi}\norm{\vphi}_{\sharp}\\&
\leq C(1+\norm{\u})\norm{\nablag\vphi}\norm{\vphi}_{\sharp}+C\norm{\nablag\vphi}^\frac12\norm{\vphi}_{\sharp}^\frac32\\&
\leq \frac 14 \norm{\nablag\vphi}^2+C(1+\norm{\u}^2)\norm{\vphi}^2_\sharp.
 \end{align*}
 To conclude, we have, by the regularity of the flow map, 
 \begin{align*}
    & \norma{\intg \vphi\mathcal{N}(t)\vphi Hv_\n\ds-\frac12\intg\norma{\nablag\mathcal{N}(t)\vphi}^2 Hv_\n\ds+\frac12\intg \nablag\mathcal{N}(t)\vphi\cdot \E_S(\VVn)\nablag\mathcal{N}(t)\vphi}\ds\\&
    \leq C(\norm{\vphi}\norm{\NN(t)\vphi}+\norm{\vphi}_\sharp^2)
    \leq C(\norm{\nablag\vphi}^\frac12\norm{\vphi}_{\sharp}^\frac32+\norm{\vphi}_\sharp^2)\\&
    \leq \frac14 \norm{\nablag\vphi}^2+C\norm{\vphi}_{\sharp}^2.
 \end{align*}
To sum up, putting all these estimates in \eqref{Gron}, we get
\begin{align}
    \frac d{dt}\norm{\vphi}_\sharp^2+\frac{1}{4}\norm{\nablag\vphi}^2+\theta_0\norm{\vphi}^2\leq C(T)(1+\norm{\u}^2)\norm{\vphi}_\sharp^2,\quad \forall t\in[0,\tT],
    \label{contest1}
\end{align}
 so that uniqueness follows by an application of the Gronwall lemma, recalling that $\u\in L^2_{\H^1(\tT)}$.
 \begin{remark}
 We point out that uniqueness holds, with the same proof, for weak solutions. Furthermore, to ensure the sole uniqueness the regularity $\u\in L^2_{\L^2(\tT)}$ is enough.
 \label{wekasolutions}
 \end{remark}
   \textbf{Existence.} 
    Very similar results have been obtained also in \cite[Theorem 3.9]{CEGP} under even more general assumptions on the normal evolution velocity $\V_\n$ of the surface. Here the main difference is related to the advective term $\u$, corresponding to set $\V_\tau=\V_a=\u$ in \cite{CEGP}, since in this case it is much less regular, compared to the assumptions on the advective terms $\mathbf{V}_a$ and $\V_\tau$ (see also assumption \cite[(2.4)]{CEGP}). Nevertheless, the approximation procedure by means of a suitable Galerkin scheme and an approximation of the singular potential $\Psi$ can be borrowed identically also in this case. In particular, since $\V_\n$ is already very regular due to the flow map, we can approximate the velocity $\u$ with a sequence $\{\u_k\}_k\subset {C}^\infty_{\H^1\cap \L^2}$. Another possibility is to show the existence of a local regular solution in a way similar to what is done in \cite[Theorem 4.1]{AGP1} and use that solutions for the computations.
    
    In any case, for the sake of brevity, here we only present the formal necessary estimates, and we refer to the proof of \cite[Theorem 3.9]{CEGP} for the more technical details on a possible rigorous approximating scheme.

    Let us then operate formally, assuming to have a sufficiently regular solution. Since in general in the approximation schemes we do not have the bound $\norma{\vphi}<1$, we do not assume this and will deduce it as a consequence of the estimates.
Now, first, we observe that the total mass is conserved. Indeed, integrating by parts, recalling \eqref{dt1},
$$
\frac{d}{dt}\intg\vphi\ds =\intg \dt\vphi\ds+\intg\vphi\divg\VVn\ds=\intg \dt\vphi\ds+\intg\vphi H\v_n\ds,
$$
and, thus by integrating over $\gam(t)$ equation \eqref{pr}$_1$, we get, integrating by parts and recalling $\divg\v=0$ and $\v=\u+\VVn$,
\begin{align*}
0=\intg \dt\vphi\ds-\intg \vphi \divg\v \ds+\intg \vphi H\v\cdot\n\ds=\intg \dt\vphi\ds+\intg \vphi H\v_n\ds=\frac d{dt}\intg\vphi\ds,
\end{align*}
entailing the desired mass conservation
\begin{align}
  \overline{\vphi}(t)=\frac{ \intg \vphi\ds}{\norma{\Gamma(t)}}=\frac{ \intg \vphi\ds}{\norma{\Gamma_0}}\equiv\frac{ \int_{\gz} \vphi_0\ds}{\norma{\Gamma_0}},\quad \forall t\in[0,T],
  \label{masscons}
\end{align}
where we have used the surface inextensibiliy assumption, so that $\norma{\Gamma(t)}=\norma{\gz}$ for any $t\in[0,T]$.
    
    Now, recalling \eqref{dt1} and \eqref{dt3}, we have, since $\divg\V_\n=H v_\n $, 
    \begin{align*}
        \frac12\frac{d}{dt}\Vert \nablag\vphi\Vert^2&=\intg \nablag \vphi \cdot \nablag\dt\vphi \ds+\frac12\intg \vert \nablag\vphi\vert^2H v_\n\ds -\intg \nablag\vphi\cdot \E_S(\V_\n)\nablag\vphi\ds,
    \end{align*}
    as well as
    \begin{align*}
        \frac d{dt} \intg \Psi(\vphi)\ds=\intg\Psi'(\vphi)\dt\vphi \ds+\intg\Psi(\vphi)H v_\n\ds .
    \end{align*}
    Therefore, we obtain a standard energy estimate by multiplying \eqref{pr}$_1$ by $\mu$: since, integrating by parts,
    $$
    \intg \dt\vphi\mu\ds=\intg \nablag\vphi\cdot \nablag\dt\vphi\ds+\intg\Psi'(\vphi)\dt\vphi\ds,
    $$
we infer
    \begin{align}
        \nonumber&\frac{d}{dt}\left(\frac12\Vert \nablag\vphi\Vert^2+\intg\Psi(\vphi)\right)\ds+\intg\vert \nablag\mu\vert^2\ds+\intg \nablag\vphi\cdot \v \ \mu\ds\\&=\frac12\intg \vert \nablag\vphi\vert^2H v_\n\ds -\intg \nablag\vphi \cdot\E_S(\V_\n)\nablag\vphi\ds+\intg\Psi(\vphi)H v_\n\ds .
        \label{dtfin}
    \end{align}
Moreover, by multiplying equation \eqref{pr}$_1$ by $\varphi$, we deduce, again using \eqref{dt1} and integrating by parts,
\begin{align*}
    \frac12\frac d {dt}\norm{\vphi}^2+\norm{\Delta_\Gamma\vphi}^2+\intg \Psi''(\vphi)\vert \nablag\vphi\vert^2\ds+\intg \vert\vphi\vert^2 H v_\n\ds =0,
\end{align*}
where we used the fact that $\divg\v=0$. This entails, by the regularity of the flow map, that
$$
\frac d{dt}\norm{\vphi}^2+\norm{\Delta_\Gamma\vphi}^2+\intg F''(\vphi)\vert \nablag\vphi\vert^2\ds\leq C\norm{\vphi}^2+\theta_0\norm{\nabla\vphi}^2,
$$
but the standard interpolation
$$
\theta_0\norm{\nablag\vphi}^2=-\theta_0\intg \vphi\Delta_{\gam}\vphi\ds\leq \theta_0\norm{\vphi}\norm{\Delta_{\gam}\vphi}
$$
entails by Young's inequality that 
$$
\frac d{dt}\norm{\vphi}^2+\frac12\norm{\Delta_\Gamma\vphi}^2+\intg F''(\vphi)\vert \nablag\vphi\vert^2\ds\leq C\norm{\vphi}^2,\quad \forall t\in[0,\tT]
$$
so that, by Gronwall's Lemma we deduce 
\begin{align}
\vphi\in L^\infty_{L^2(\tT)}.
    \label{vphiinft}
\end{align}
Now we can estimate all the terms appearing in the identity \eqref{dtfin}. First, integrating by parts, by Poincaré's and Sobolev-Gagliardo-Nirenberg's inequalities, recalling the regularity of $\VVn$, \eqref{masscons} and \eqref{vphiinft}, we deduce that 
\begin{align*}
    \norma{\intg \nablag\vphi\cdot \v \ \mu\ds}&\leq \norma{\intg \vphi \v\cdot \nablag \mu\ds}+\norma{\intg\vphi\mu Hv_\n\ds}\\&\leq \norm{\vphi}_{L^4\gt}(1+\norm{\u}_{\L^4\gt})\norm{\nablag\mu}+C\norm{\vphi}\norm{\mu} \\&\leq \norm{\vphi}^\frac12\norm{\vphi}_{H^1(\gam(t))}^\frac12(1+\norm{\u}_{\L^4\gt})\norm{\nablag\mu}+C\norm{\mu}\\&
    \leq C(1+\norm{\nablag\vphi}^\frac12)(1+\norm{\u}_{\L^4\gt})\norm{\nablag\mu}+C(\norm{\mu-\overline{\mu}}+\norma{\overline{\mu}})\\&
    \leq C_\epsilon(1+\norm{\u}_{\L^4\gt}^2)(1+\norm{\nablag\vphi}^2)+\frac 12 \norm{\nablag\mu}^2+\epsilon\norma{\overline{\mu}}^2,
\end{align*}
for some $\epsilon>0$ to be chosen later on. Concerning the other terms, we easily get, by the regularity of the flow map, 
\begin{align*}
    &\norma{\frac12\intg \vert \nablag\vphi\vert^2H v_\n\ds -\intg \nablag\vphi \cdot\E_S(\V_\n)\nablag\vphi\ds+\intg\Psi(\vphi)H v_\n }\ds\\&\leq C\left(\frac12\Vert \nablag\vphi\Vert^2+\intg(\Psi(\vphi)
+\tilde{C})\ds\right)+C,
\end{align*}
    where $\tilde{C}>0$ is a suitable constant such that $\Psi\geq -\tilde{C}$ (recall that $\Psi$ is bounded below by assumption, and can be chosen as such also in a suitable approximation). We are then left with the estimation of $\overline{\mu}$, which can be performed similarly to \cite{DE,CEGP}. In particular, thanks to mass conservation, it holds
    \begin{equation}
				\| F^{\prime }(\vphi )\| _{L^{1}(\gam(t) )}\leq
				C_1\left\vert \intg F^{\prime }(\vphi )\big(%
				\vphi -\overline{\vphi }\big)\ds\right\vert
				+C_{2},  \label{MZ-2}
			\end{equation}%
			where {$C_1$ and $C_2$ } are positive constants only depending on $F'$, $\norma{\Gamma_0}$ and $\overline{\vphi}_{0}$. We sketch the argument (devised, for instance, in \cite{Kenmochi}) to show that the constants are not depending on other quantities. Since by assumption it holds $\norma{\overline{\vphi}}\in(0,1)$, and since $F'(s)\to \pm\infty$ as $s\to \pm1$ and $F'$ is monotone increasing, there exists $m\in[0,1)$ such that $F'(s)>0$ for any $s\in [m,1)$, $F'(s)<0$ for any $s\in (-1,-m]$ and $\norma{\overline{\vphi}_0}<m$. Then we set
   $$
   \Gamma_1(t):=\{x\in \gam(t):\ -m\leq \vphi(x,t)\leq m\},\quad 
   \Gamma_2(t):=\{x\in \gam(t):\ \vphi(x,t)< -m \text{ or }\vphi(x,t)>m\}, 
   $$
and define $\delta_1:=\min\{\overline{\vphi}_0+m,m-\overline{\vphi}_0\}>0$. Then it holds, recalling that $\norma{\gam(t)}=\norma{\gz}$ and $\overline{\vphi}=\overline{\vphi}_0$,
\begin{align*}
    \delta_1\norm{F'(\vphi)}_{L^1(\Gamma(t))}&= \delta_1\ints{{\Gamma_1}}\norma{F'(\vphi)}\ds+\delta_1\ints{{\Gamma_2}}\norma{F'(\vphi)}\ds\\&
\leq \delta_1\max_{s\in[-m,m]}\norma{F'(s)}\norma{\gam(t)}+\ints{\gam_2}{F'(\vphi)}(\vphi-\overline{\vphi}_0)\ds\\&
\leq 
\delta_1\max_{s\in[-m,m]}\norma{F'(s)}\norma{\gam(t)}+\ints{\gam}{F'(\vphi)}(\vphi-\overline{\vphi}_0)\ds-\ints{\gam_1}{F'(\vphi)}(\vphi-\overline{\vphi}_0)\ds
\\&\leq 2\delta_1\max_{s\in[-m,m]}\norma{F'(s)}\norma{\gz}+\norma{\ints{\gam}{F'(\vphi)}(\vphi-\overline{\vphi})\ds},
\end{align*}
where in the last inequality we used the fact that 
$$
-\ints{\gam_1}{F'(\vphi)}(\vphi-\overline{\vphi}_0)\ds\leq \delta_1\max_{s\in[-m,m]}\norma{F'(s)}\norma{\gz}.
$$ Therefore, it is easy to identify $C_1:=\frac1{\delta_1}$ and $C_2:=2\max_{s\in[-m,m]}\norma{F'(s)}\norma{\gz}$ and these constants only depend on $\overline{\vphi}_0, F', \norma{\gz}$ as expected. Now, from the definition of $\mu$, integrating by parts, recalling \eqref{vphiinft}, we get
\begin{align}
\norm{\nablag\vphi}^2+\ints{\gam}{F'(\vphi)}(\vphi-\overline{\vphi})\ds\leq \norm{\mu-\overline{\mu}}\norm{\vphi-\overline{\vphi}}+\theta_0\norm{\vphi}\norm{\vphi-\overline{\vphi}}\leq C(1+\norm{\mu-\overline{\mu}}),
    \label{mutest}
\end{align}
and thus, by Young's and Poincaré's inequalities, and exploiting \eqref{MZ-2}, we deduce 
\begin{align}
   \norma{\overline{   \mu}}=\frac{1}{\norma{\gz}} \| F^{\prime }(\vphi )\| _{L^{1}(\gam(t) )}\leq
				\frac{1}{\norma{\gz}}C_1\left(C_3\norm{\nablag\mu}+C\right)
				+\frac{1}{\norma{\gz}}C_{2},
    \label{mubar}
\end{align}
for some $C_3>0$, and thus 
\begin{align}
   \norma{\overline{   \mu}}^2\leq
				\frac{2}{\norma{\gz}^2}C_1^2C_3^2\norm{\nablag\mu}^2+\frac{2}{\norma{\gz}^2}\left(C_1C
				+C_{2}\right)^2.
    \label{muep2}
\end{align}
Therefore, from \eqref{dtfin} and all the inequalities above, by choosing $\epsilon=\frac{1}{8\norma{\gz}^2(C_1C_3)^2}$we end up with 
\begin{align*}
   & \frac{d}{dt}\left(\frac12\Vert \nablag\vphi\Vert^2+\intg(\Psi(\vphi)+\tilde{C})\ds\right)+\frac14\intg\vert \nablag\mu\vert^2\ds\\&\leq C(1+\norm{\u}_{\L^4(\gam(t))}^2)\left(\frac12\Vert \nablag\vphi\Vert^2+\intg(\Psi(\vphi)+\tilde{C})\ds\right)+C(1+\norm{\u}_{\L^4(\gam(t))}^2),\quad \forall t\in[0,\tT].
\end{align*}
Therefore, an application of Gronwall's Lemma and \eqref{mubar}, assuming that $\norm{\u}_{L^2_{\H^1(\tT)}}\leq C_A$, for some $C_A>0$, entail 
\begin{align}
    &\norm{\vphi}_{L^\infty_{H^1(\tT)}}+
    \norm{\mu}_{ L^2_{H^1(\tT)}}+\norm{\Psi(\vphi)}_{L^\infty_{L^1(\tT)}}\leq C(T,C_A)
    \label{ppsi}.
\end{align}
Observe now that, by testing the equation for $\mu$ by $-\Delta_{\gam(t)}\vphi$ we immediately deduce, by standard estimates, since $F''\geq\theta$, 
\begin{align*}
\norm{\Delta_{\gam(t)}\vphi}^2\leq \norm{\Delta_{\gam(t)}\vphi}^2+\intg F''(\vphi)\norma{\nablag\vphi}^2\ds\leq C(\norm{\nablag\mu}\norm{\nablag\vphi}+\norm{\nablag\vphi}^2),
\end{align*}
so that we also deduce
$$
\norm{\vphi}_{L^4_{H^2(\tT)}}\leq C(T,C_A),
$$
from which we also deduce by comparison in the equation for $\mu$ that 
\begin{align}
\norm{F'(\vphi)}_{L^2_{L^2(\tT)}}\leq C(T,C_A).
\label{fpp}
\end{align}
This last regularity result also allows to infer in the end, by standard arguments (after passing to the limit in the complete approximating scheme), that $\vert \vphi\vert<1$ almost everywhere on $\Gamma(t)$, for any $t\in[0,T]$. 

In conclusion, by comparison it is now immediate to deduce, recalling $\norma{\vphi}<1$, that 
$$
\norm{\dt\vphi}_{H^1(\gam(t))'}\leq C(\norm{\nablag\mu}+\norm{\v}),
$$
and thus also 
$$
\norm{\vphi}_{H^1_{H^{-1}(\tT)}}\leq C(T,C_A).
$$
The regularity shown so far is enough to show the existence of a weak solution (notice that the initial data can be less regular, in particular it is enough that $\vphi_0\in H^1(\gz)$, $\norma{\overline{\vphi}_0}<1$ and $\norm{\vphi_0}_{L^\infty(\gz)}\leq1$). We now aim at showing higher-order estimates (always with the idea of being in the same Galerkin framework with the use of a potential $\Psi$ approximating the singular one as in \cite{CEGP}). Let us multiply equation \eqref{pr}$_1$ by $\dt\mu$. Again from \eqref{dt3} we infer
\begin{align*}
    \frac12\frac d{dt}\norm{\nablag\mu}^2=\intg \nablag\dt\mu\cdot\nablag\mu\ds+\intg\norma{\nablag\mu}^2Hv_\n\ds-2\intg\nablag\mu\cdot\E_s(\V_\n)\nablag\mu\ds.
\end{align*}
Therefore we have
\begin{align}
    \nonumber&\frac12 \ddt\norm{\nablag\mu}^2+\intg \dt\vphi \dt\mu\ds+\intg \nablag\vphi\cdot \v\dt\mu\ds\\&=\intg\norma{\nablag\mu}^2Hv_\n\ds-2\intg\nablag\mu\cdot\E_s(\V_\n)\nablag^T\mu\ds.\label{dmu}
\end{align}
Now, we can argue in the same way as to obtain \cite[(2.36)]{CEGP}, namely, assuming $\eta\in C^\infty_{H^1}$, we can write, exploiting \eqref{dt1} and \eqref{dt3},
\begin{align}
    &\frac{d}{dt}\intg\nablag\vphi\cdot\nablag\eta\ds+\frac{d}{dt}\intg\Psi'(\vphi)\eta\ds=\frac{d}{dt}\intg \mu\eta\ds\non\\&=\intg \dt\mu\eta\ds+\intg\mu\dt\eta\ds+\intg\mu\eta\divg\VVn\ds,\label{mmm1}
\end{align}
but, recalling the definition of $\mu$, 
\begin{align*}
&\frac{d}{dt}\intg\nablag\vphi\cdot\nablag\eta\ds+\frac{d}{dt}\intg\Psi'(\vphi)\eta\ds\\&=\intg \nablag \dt\vphi\cdot \nablag\eta\ds+\intg \Psi''(\vphi)\dt\vphi\eta\ds+\intg \mu\dt\eta\ds\\&+\intg \nablag\vphi\cdot\nablag \eta \divg\VVn\ds+\intg\Psi'(\vphi)\eta\divg\VVn\ds-2\intg \nablag\vphi\cdot\E_S(\VVn)\nablag\eta\ds,
\end{align*}
so that, inserting in \eqref{mmm1}, we deduce
\begin{align}
   \nonumber \intg\dt\mu\eta\ds&=-\intg\mu\eta\divg\VVn\ds+\intg \nablag \dt\vphi\cdot \nablag\eta\ds+\intg \Psi''(\vphi)\dt\vphi\eta\ds\\&+\intg \nablag\vphi\cdot\nablag \eta \divg\VVn\ds+\intg\Psi'(\vphi)\eta\divg\VVn\ds-2\intg \nablag\vphi\cdot\E_S(\VVn)\nablag\eta\ds,\label{muesssential}
\end{align}
and then by density it is immediate to deduce that this holds for any $\eta\in L^2_{H^1}$. By setting $\eta=\dt\vphi$ we can then estimate the term $\intg\dt\vphi\dt\mu$. In the end we thus obtain from \eqref{dmu}, recalling the definition of $\Psi$ and that $\divg\VVn=H v_\n$, that
\begin{align}
    \nonumber&\frac12 \frac d{dt}\norm{\nablag\mu}^2+\norm{\nablag\dt\vphi}^2+
    (\dt\vphi,F''(\vphi)\dt\vphi)+\intg \nablag\vphi\cdot \v\dt\mu\ds\\&=\theta_0\norm{\dt\vphi}^2+\intg\norma{\nablag\mu}^2Hv_\n\ds-2\intg\nablag\mu\cdot\E_s(\V_\n)\nablag\mu\ds  \label{mmmua}\\&\non
    +\intg \mu\dt\vphi Hv_\n\ds-\intg\nablag\vphi\cdot\nablag\dt\vphi Hv_\n\ds-\intg \Psi'(\vphi)\dt\vphi Hv_\n\ds+2\intg \nablag\vphi\cdot\E_S(\VVn)\nablag\dt\vphi\ds.
\end{align}
Now we observe that, by the regularity of the flow map, and estimate \eqref{mubar}, we have
\begin{align*}
    &\norma{\intg\norma{\nablag\mu}^2Hv_\n\ds-2\intg\nablag\mu\cdot\E_s(\V_\n)\nablag\mu
    \ds+\intg \mu\dt\vphi Hv_\n\ds}\leq C(1+\norm{\nablag\mu}^2)+C\norm{\dt\vphi}^2,
\end{align*}
and, similarly, recalling that $\vphi\in L^\infty_{H^1}$,
\begin{align*}
&\norma{-\intg\nablag\vphi\cdot\nablag\dt\vphi Hv_\n\ds+2\intg \nablag\vphi\cdot\E_S(\VVn)\nablag\dt\vphi\ds}\leq C+\frac12\norm{\nabla\dt\vphi}^2.
\end{align*}
Now we need to control the term with $\Psi'$. We follow \cite{CEGP}, showing some estimates also working in the Galerkin setting. In particular, observe that, by \eqref{dt1},
\begin{align*}
    \frac{d}{dt}\intg \Psi(\vphi)Hv_\n\ds=\intg \Psi'(\vphi)\dt\vphi Hv_\n\ds+\intg \Psi(\vphi)\dt(Hv_\n)\ds+\intg \Psi(\vphi)(Hv_\n)^2\ds,
\end{align*}
and also, by the regularity of the flow map, we have
\begin{align*}
    \norma{\intg \Psi(\vphi)\dt(Hv_\n)\ds+\intg \Psi(\vphi)(Hv_\n)^2}\ds\leq C\norm{\Psi(\vphi)}_{L^1(\Gamma(t))}\leq C,
\end{align*}
recalling \eqref{ppsi}. Therefore, recalling that $F''\geq \theta$, \eqref{mmmua} becomes
\begin{align}
    \nonumber& \frac d {dt}\left(\frac12\norm{\nablag\mu}^2+\intg \Psi(\vphi)Hv_\n\right)\ds+\frac12
\norm{\nablag\dt\vphi}^2+
    \theta\norm{\dt\vphi}^2+\intg \nablag\vphi\cdot \v\dt\mu\ds\\&\leq \tilde{C}\norm{\dt\vphi}^2+C(1+\norm{\nablag\mu}^2),
    \label{mmmuaa}
\end{align}
for some $\tilde{C}>0$.
In order to control $\norm{\dt\vphi}^2$, we test \eqref{pr}$_1$ by $\dt\vphi$ and integrate by parts. This gives, by standard estimates, recalling $\vphi\in L^\infty_{H^1}$,
\begin{align}
    \tilde{C}\norm{\dt\vphi}^2&\leq \norm{\nablag\mu}\norm{\nablag\dt\vphi}+\norm{\vphi}_{L^4(\gam(t))}\norm{\v}_{\L^4(\gam(t))}\norm{\nablag\dt\vphi}\nonumber\\&
    \leq \frac{1}{8}\norm{\nablag\dt\vphi}^2+C(\norm{\nablag\mu}^2+\norm{\v}_{\L^4(\gam(t))}^2),\label{ddtt}
\end{align}

Now, to estimate the advective term related to $\v$, we have to work a bit more. First, if we consider the Galerkin setting, together with the approximation $\u_k\in C^\infty_{\H^1\cap \L^2}$ of the velocity (here still denoted for simplicity by $\v$), we have, by \eqref{dt1}-\eqref{dt2},
\begin{align*}
   & {\intg \nablag\vphi\cdot \v\dt\mu\ds}=\frac{d}{dt}\intg \nablag\vphi\cdot \v\mu\ds-\intg \nablag\dt\vphi\cdot \v\mu\ds-\intg \nablag\vphi\cdot\dt \v\mu\ds\\&-\intg \nablag\vphi\cdot \v\mu Hv_\n\ds+\intg \v\mu\cdot (\nablag^T\VVn)\nablag\vphi\ds,
\end{align*}
and by standard estimates, recalling \eqref{muep2} and the regularity of the flow map, we deduce
\begin{align*}
    &\norma{\intg \nablag\dt\vphi\cdot \v\mu\ds+\intg \nablag\vphi\cdot\dt \v\mu\ds+\intg \nablag\vphi\cdot \v\mu Hv_\n\ds-\intg \v\mu\cdot (\nablag^T\VVn)\nablag\vphi\ds}\\&\leq  C(1+\norm{\dt\v}_{\L^\infty(\gam(t))}+\norm{\v}_{\L^\infty(\gam(t))}^2)(1+\norm{\nablag\mu}^2)+\frac18 \norm{\nablag\dt\vphi}^2,
\end{align*}
and thus we obtain from \eqref{mmmuaa} that
\begin{align}
    \nonumber& \frac d {dt}\left(\frac12\norm{\nablag\mu}^2+\intg \Psi(\vphi)Hv_\n\ds+\intg \nablag\vphi\cdot \v\mu\ds\right)+\frac14
\norm{\nablag\dt\vphi}^2+
    \theta\norm{\dt\vphi}^2\\&\leq C(1+\norm{\dt\v}_{\L^\infty(\gam(t))}+\norm{\v}_{\L^\infty(\gam(t))}^2)(1+\norm{\nablag\mu}^2).
    \label{mmmua2}
\end{align}
Since $\divg\v=0$, it holds
$$
\norma{\intg \nablag\vphi\cdot \v\mu\ds}=\norma{-\intg \vphi \v\cdot\nablag\mu\ds}\leq \frac{1}{4}\norm{\nablag\mu}^2-C\norm{\v}_{\L^\infty(\gam(t))}^2,
$$
and also, by \eqref{ppsi}, we have
$$
\norma{\intg \Psi(\vphi)Hv_\n\ds}\leq C.
$$
Then, an application of Gronwall's lemma, together with \eqref{mubar}, ensures that $\mu\in L^\infty_{H^1(\tT)}$ and $\vphi\in H^1_{H^1(\tT)}$. As a matter of facts, this result depends on the approximating sequence $\u_k$ and the approximation of the singular potential $\Psi$, but it is independent of the Galerkin parameters. Therefore, we can easily pass to the limit in these Galerkin parameters and think of finding uniform estimates in the parameter $k$ of the approximating velocities and also uniform with respect to the approximation of $\Psi$. 

We now start again from \eqref{mmmua}, assuming a solution after passing to the limit in the Galerkin scheme, which is sufficiently regular to perform the estimate, since the singular potential is still approximated by a regular one and the velocity is still the approximated one (i.e., $\u_k$). We are allowed to perform estimates in the $L^p$ setting, which are not allowed in a Galerkin scheme. In particular, exactly with the same arguments as in \cite[(3.60)]{CEGP}, by testing the equation for $\mu$ by $F'(\vphi)\norma{F'(\vphi)}^{p-2}$, for $p\geq 2$ and applying suitable H\"{o}lder and Young's inequalities, we can deduce, also recalling \eqref{mubar}, that 
\begin{align}
    \nonumber&\norm{F'(\vphi)}_{L^p(\gam(t))}^p+\intg(p-1)F''(\vphi)\norma{F'(\vphi)}^{p-2}\norma{\nabla\vphi}^2\ds\\&\non\leq \norm{\mu+\theta_0\vphi}_{L^p(\gam(t))}\norm{F'(\vphi)}^{p-1}_{L^p(\gam(t))}\\&
    \leq \frac 12\norm{F'(\vphi)}_{L^p(\gam(t))}^p+C(p)\norm{\mu+\theta_0\vphi}_{L^p(\gam(t))}^p,\label{Fprimeest}
\end{align}
entailing, since $H^1(\gam(t))\hookrightarrow L^p(\gam(t))$ for any $p\geq2$ ad recalling \eqref{mubar}, that
\begin{align}
\norm{F'(\vphi)}_{L^p(\Gamma(t))}\leq C(p)(1+\norm{\nablag\mu}),
\label{F'}
\end{align}
which also entails, since $\vphi\in L^\infty_{H^1(\tT)}$, 
\begin{align}
\norm{\Psi'(\vphi)}_{L^p(\Gamma(t))}\leq C(p)(1+\norm{\nablag\mu}).
\label{Psi'}
\end{align}
From the equation for $\mu$ we thus obtain, by elliptic regularity,
\begin{align}
\norm{\vphi}_{W^{2,p}(\gam(t))}\leq  C(p)(1+\norm{\nablag\mu}), \quad \forall \in [2,+\infty).
    \label{W2p}
\end{align}
Now we can analyze the remaining advective term \eqref{mmmuaa}. Thanks to \eqref{muesssential}, with $\eta=\nablag\vphi\cdot \v$, we infer
\begin{align}
   \nonumber &\intg\dt\mu\nablag\vphi\cdot \v\ds=-\intg\mu\nablag\vphi\cdot \v\divg\VVn\ds+\intg \nablag \dt\vphi\cdot \nablag(\nablag\vphi\cdot \v))\ds\\&+\intg \Psi''(\vphi)\dt\vphi\nablag\vphi\cdot \v\ds\nonumber+\intg \nablag\vphi\cdot\nablag (\nablag\vphi\cdot \v) \divg\VVn\ds\\&+\intg\Psi'(\vphi)\nablag\vphi\cdot \v\divg\VVn\ds-2\intg \nablag\vphi\cdot\E_S(\VVn)\nablag(\nablag\vphi\cdot \v)\ds.\label{muesssential2}
\end{align}
We can now estimate all the terms appearing in the right-hand side. First, we have, recalling \eqref{mubar} and \eqref{W2p} with $p=2$, and the regularity of the flow map,
\begin{align*}
    &\norma{\intg\mu\nablag\vphi\cdot \v\divg\VVn\ds}\leq C\norm{\mu}\norm{\nablag\vphi}_{\L^4(\gam(t))}\norm{\v}_{\L^4(\gam(t))}\\&\leq C(1+\norm{\nablag\mu})\norm{\vphi}_{H^2(\gam(t))}(1+\norm{\u}_{\L^4(\gam(t))})\leq C(1+\norm{\nablag\mu}^2)(1+\norm{\u}_{\L^4(\gam(t))}).
\end{align*}
Then, we have, recalling the regularity of the flow map and $\v=\u+\VVn$, the embedding $W^{2,4}\gt\hookrightarrow W^{1,\infty}\gt$ and \eqref{W2p} with $p=4$, 
\begin{align*}
    &\norma{\intg \nablag \dt\vphi\cdot \nablag(\nablag\vphi\cdot \v))\ds}+\norma{\intg \nablag\vphi\cdot\nablag (\nablag\vphi\cdot \v) \divg\VVn\ds}\\&+\norma{2\intg \nablag\vphi\cdot\E_S(\VVn)\nablag(\nablag\vphi\cdot \v)\ds}\\&\leq C(\norm{\nablag\dt\vphi}+\norm{\nablag\vphi})\left(\norm{\nablag\vphi}_{\L^\infty(\gam(t))}(1+\norm{\u}_{\H^1\gt})+\norm{\vphi}_{W^{2,4}(\gam(t))}(1+\norm{\u}_{\L^4(\gam(t))})\right)\\&
    \leq C(1+\norm{\nablag\dt\vphi})\left(\norm{\vphi}_{W^{2,4}(\gam(t))}(1+\norm{\u}_{\H^1\gt})+\norm{\vphi}_{W^{2,4}(\gam(t))}(1+\norm{\u}_{\L^4(\gam(t))})\right)\\&
    \leq C(1+\norm{\nablag\dt\vphi})(1+\norm{\nablag\mu})(1+\norm{\u}_{\H^1\gt}+\norm{\u}_{\L^4(\gam(t)})\\&
    \leq 
    \frac{1}{16}\norm{\nablag\dt\vphi}^2+\norm{\nabla\mu}^2(1+\norm{\u}_{\H^1\gt}^2+\norm{\u}_{\L^4(\gam(t))}^2)+C(1+\norm{\u}_{\H^1\gt}^2+\norm{\u}_{\L^4(\gam(t))}^2).
\end{align*}
Then we have, integrating by parts recalling $\divg\v=0$ and \eqref{Psi'}, and by H\"{o}lder's and Young's inequalities,
\begin{align*}
    \norma{\intg \Psi''(\vphi)\dt\vphi\nablag\vphi\cdot \v\ds}&= \norma{\intg \nablag\Psi'(\vphi)\dt\vphi\cdot \v\ds}\\&\leq\norma{\intg \Psi'(\vphi)\nablag\dt\vphi\cdot \v\ds}+\norma{\intg \Psi'(\vphi)\dt\vphi Hv_\n\ds}\\&
    \leq \norm{\nablag\dt\vphi}\norm{\Psi'(\vphi)}_{L^4\gt}\norm{\v}_{\L^4\gt}+C\norm{\Psi'(\vphi)}\norm{\dt\vphi}\\&
\leq \frac{1}{16}\norm{\nablag\dt\vphi}^2+\frac\theta2\norm{\dt\vphi}^2+C(1+\norm{\nablag\mu}^2)(1+\norm{\u}_{\L^4\gt}^2).
\end{align*}
Proceeding in the estimates, we are left with 
\begin{align*}
\norma{\intg\Psi'(\vphi)\nablag\vphi\cdot \v\divg\VVn\ds}&\leq C\norm{\Psi'(\vphi)}_{L^4\gt}\norm{\nablag\vphi}\norm{\v}_{\L^4\gt}\\&\leq 
C\norm{\Psi'(\vphi)}_{L^4\gt}\norm{\v}_{\L^4\gt}\\&\leq C\norm{\nablag\mu}\norm{\v}_{\L^4\gt}+C(1+\norm{\u}_{\L^4\gt}),
\end{align*}
where again we exploited \eqref{Psi'}. To sum up all these estimates, coming back to \eqref{mmmuaa} and using also \eqref{ddtt}, we end up with 
\begin{align}
    \nonumber& \frac d {dt}\frac12\norm{\nablag\mu}^2+\frac14
\norm{\nablag\dt\vphi}^2+
    \frac\theta2\norm{\dt\vphi}^2\\&\leq C\norm{\nablag\mu}^2(1+\norm{\u}^2_{\H^1\gt}+\norm{\u}_{\L^4\gt}^2)+C(1+\norm{\u}^2_{\H^1\gt}+\norm{\u}_{\L^4\gt}^2)
    \label{mmmuaab}
\end{align}
from which, since $\u\in L^2_{\H^1(\tT)}$, $\norm{\u}_{L^2_{\H^1(\tT)}}\leq C_A$, and $\H^1\gt\hookrightarrow \L^4\gt$, by Gronwall's Lemma we infer, since $\mu_0\in H^1(\gz)$ and recalling once more \eqref{muep2}, that 
\begin{align}
\norm{\mu}_{L^\infty_{H^1(\tT)}}+\norm{\dt\vphi}_{L^2_{H^1(\tT)}}\leq C(T,C_A),
    \label{newreg}
\end{align}
this time independently on \textit{any} approximating parameter.
Also, thanks to \eqref{F'} and \eqref{W2p}, we get, for any $p\in[2,\infty)$,
\begin{align}
\norm{\vphi}_{L^\infty_{W^{2,p}(\tT)}}+\norm{F'(\vphi)}_{L^\infty_{L^p(\tT)}}\leq C(p,T,C_A).
    \label{2p}
\end{align} 
 Now, by elliptic regularity (see also \eqref{uniformelliptic}) applied to \eqref{pr}$_1$ we deduce, recalling $\vphi\in L^\infty_{W^{2,p}(\tT)}$ for any $p\geq2$, and $\mu \in L^\infty_{H^1(\tT)}$,
\begin{align*}
 \norm{\mu}_{H^3(\Gamma(t))}\leq C(T)({1}+\norm{\dt\vphi}_{H^1\gt}+\norm{\v\cdot\nabla_\Gamma\vphi}_{H^1\gt})\leq C(T)({1}+\norm{\dt\vphi}_{H^1\gt}+\norm{\u}_{\H^1\gt}),
\end{align*}
so that $\mu\in L^2_{H^3(\tT)}$.
All these estimates allow to conclude the regularity part of the theorem, since now we can pass to the limit in the approximating sequence $\{\u_k\}_k$ and in the approximating sequence of regular potentials to conclude. We are left to discuss the strict separation property, which can be retrieved \textit{a posteriori}. In particular, by repeating the computations leading to \eqref{Fprimeest} (this time for the solution without any approximating scheme) and, recalling that, since we are using the singular potential $\Psi$, it holds $\norma{\vphi}<1$ for almost any $t$ and $x\in\Gamma(t)$, one more precisely obtains, recalling \eqref{Gagliardo} and \eqref{newreg},
\begin{align}
\nonumber\norm{\mu}_{L^p\gt}+\norm{F'(\vphi)}_{L^p\gt}&\leq C(1+\norm{\mu}_{L^p\gt})\\&\nonumber\leq C(T)\sqrt{p}(1+\norm{\mu}_{H^1\gt})\\&\leq C(T,C_A)\sqrt p\quad \forall p\in[2,+\infty),\quad \text{ for a.a. }t\in[0,\tT].\label{ppp} 
\end{align}
Thanks to this inequality and assumption \eqref{assmp} on $F'$, and recalling \eqref{Gagliardo}, we can reproduce word by word the De Giorgi iteration scheme applied on the equation for the chemical potential $$\mu(t)=-\Delta_\Gamma\vphi(t)+\Psi(\vphi(t)),$$ for almost any $t\in[0,\tT]$ as in the proof of \cite[Theorem 3.3]{GalP}. Note that the scheme can be performed \textit{at any} fixed time $t\geq0$ for which the uniform estimate \eqref{ppp} is valid, since it is based only on the definition of $\mu(t)$ and on the Gagliardo-Nirenberg's inequality \eqref{Gagliardo} with uniform constants on $[0,T]$ (see inequality \cite[(3.20)]{GalP} for its application). Therefore, since we have assumed $\mu_0\in H^1(\gz)$, \eqref{Gagliardo} and \eqref{ppp} also hold for $t=0$, and thus in the proof of \cite[Theorem 3.3]{GalP} we can also obtain the  result for $t=0$, and deduce that there exists $\delta=\delta(T,C_A)>0$ such that 
\begin{align}
\sup_{t\in[0,\tT]}\norm{\vphi}_{C\gt}\leq 1-\delta.\label{sepp}
\end{align}
The last regularity $\vphi\in L^\infty_{H^3(\tT)}\cap L^2_{H^4(\tT)}$ is obtained by elliptic regularity from
$$
-\Delta_\gam\vphi=-\Psi'(\vphi)+\mu,
$$
and, thanks to the regularity of $\Psi$, 
$$
\norm{-\Psi'(\vphi)+\mu}_{L^\infty_{H^1(\tT)}}+\norm{-\Psi'(\vphi)+\mu}_{L^2_{W^{2,p}(\tT)}}\leq C(T,C_A),
$$
which holds thanks to the separation property. In conclusion, again setting this estimate in a suitable approximation scheme (or, for instance, by means of increment quotients), by standard estimates we deduce from \eqref{muesssential} that 
\begin{align*}
\norma{ <\dt\mu,\eta>_{H^{-1},H^{1}}}\leq C(\norm{\mu}+\norm{\nablag\dt\vphi}+\norm{\Psi''(\vphi)}_{L^\infty\gt}\norm{\dt\vphi}+\norm{\nablag\vphi}+\norm{\Psi'(\vphi)}_{L^\infty\gt})\norm{\eta}_{H^1\gt},
\end{align*}
for all $\eta\in H^1\gt$. Therefore, we immediately infer from the regularities above that 
$$
\norm{\dt\mu}_{L^2_{H^{-1}(\tT)}}\leq C(T,C_A),
$$
concluding the proof of the theorem.
\end{proof}
\section{Proof of Theorem \ref{strong1}}\label{sec:proof}
We divide the proof in many steps. First, we show the local existence of a strong solution according to Definition \ref{strong}. We exploit the result of \cite[Theorem 4.1]{AGP1}.

\textbf{Step 0. Existence of a local strong solution on an interval $[0,T_0]$. }
Observe that, thanks to Remark \ref{necessaryrem}, the initial data are such that 
$$
\v_0\in \H^1(\gz),\quad  \vphi_0\in H^3(\gz),$$
and 
$$\exists \delta_0\in(0,1):\ \norm{\vphi}_{L^\infty(\gz)}\leq 1-2\delta_0.
$$
Then there exist $p>4$ and $q\in(2,4)$ such that  
\begin{align}
H^3(\gz)\hookrightarrow W^{4+\epsilon-\frac 4 q,p}(\gz)\hookrightarrow B^{4-\frac 4 q}_{p,q}(\gz),
\label{H3b}
\end{align}
for $\epsilon>0$ sufficiently small (for instance,  we can choose $p=\tfrac92$, $q=\tfrac 52$ and $\epsilon=\tfrac2{45}$).
Thus the initial data and the functions $\Psi,\nu,\rho$ are sufficiently regular to apply \cite[Theorem 4.1]{AGP1}, recalling \cite[Remarks 4.3-4.4]{AGP1}. This means that there exists $0<\tT\leq T$ and $(\tilde{\v},\pi,\vphi,\mu)$, solution to \eqref{mainp}, such that 
\begin{align*}
    \tilde{\v}\in L^\infty_{\L^2}\cap L^2_{\H^2(\tT)}\cap H^1_{\L^2(\tT)},\quad \pi\in L^2_{H^1(\tT)},\quad \vphi\in L^2_{W^{4,p}(\tT)}\cap H^1_{L^p(\tT)},\quad \mu\in L^2_{W^{2,p}(\tT)},
\end{align*}
and
\begin{align}
\sup_{t\in[0,\tT]}\norm{\vphi(t)}_{L^\infty\gt}\leq 1-\delta_0.
\label{separation}
\end{align}

Now consider a strong solution $\vphi_2$ to \eqref{pr} with $\v:= \tilde{\v}+\VVn$ and initial datum $\vphi_0$. Since the velocity $\tilde{\v}$ is sufficiently regular and since, by construction, $\divg(\tilde{\v}+\VVn)=0$, this solution exists by Theorem \ref{existence}. Now we can set $\vphi_1=\vphi$ and repeat word by word the proof leading to the continuous dependence estimate \eqref{contest1}. This entails $\vphi\equiv \vphi_2$, so that we also deduce the  regularity 
$$
\vphi\in L^\infty_{W^{2,p}(\tT)}\cap H^1_{H^1(\tT)},\quad \mu\in L^\infty_{H^1(\tT)}\cap L^2_{H^3(\tT)},\quad \forall p\in[2,+\infty).
$$
Therefore we have shown that $(\tilde{\v},\pi,\vphi)$ is a local strong solution to \eqref{mainp} on $[0,T_0]$ (with $T_0=\tT$), according to Definition \ref{strong}, concluding the proof of local in time existence of a strong solution.

\textbf{Step 1. A problem for a divergence-free velocity.}
We now need to introduce a problem for a divergence-free velocity, as done in \cite[Section 5]{AGP1}. Namely, we solve the problem
\begin{align*}
    -\Delta_\Gamma \Pi=Hv_\n
\end{align*}
for $\Pi$ such that $\overline{\Pi}=0$. This admits a unique solution by elliptic theory.  We then define $\widehat{\u}=\nablag\Pi$ and set $\v=\V+\widehat{\u}$. Then, thanks to \cite[Lemma A.1]{AGP1}, the regularity on the flow map ensures that we have $\Pi\in C^1_{W^{4,p}}$, and then  $\widehat{\u}\in C^0_{\W^{3,p}}$ with $\P\dt\widehat{\u}\in C^0_{\W^{3,p}}$, for any $p\in[2,+\infty)$, meaning also that $\widehat{\u}\in C^0_{\C^2}$ as well as $\P\dt\widehat{\u}\in C^0_{C^2}$. In this way, we obtain the following equivalent problem: find $(\V,\pi,\vphi)$ such that
 \begin{align}
     \begin{cases}
     \rho \P\dt \V +((\rho\V+\J_\rho)\cdot\nabla_\Gamma)\V+\rho(\widehat{\u}\cdot\nabla_\Gamma)\V+\rho(\V\cdot\nabla_\Gamma)\widehat{\u}+\rho \Vn\mathbf{H}\V-2\P\divg(\nu(\vphi)\E_S(\V))\\+\Vn\H\J_\rho+\nablag \pi
     =-\P\divg(\nablag\vphi\otimes\nablag\vphi)+2\P\divg(\nu(\vphi)\Vn\mathbf{H})+\frac\rho2\nablag(\Vn)^2\\
    \ \quad\quad\quad\quad\quad\quad\quad\quad-\rho \P\dt \widehat{\u} -((\rho\widehat{\u}+\J_\rho)\cdot\nabla_\Gamma)\widehat{\u}-\rho \Vn\mathbf{H}\widehat{\u}+2\P\divg(\nu(\vphi)\E_S(\widehat{\u})),\\
         \divg\V=0,
         \\\dt\vphi-\Delta_{\Gamma}\mu+\nablag\vphi\cdot \V
         +\nablag\vphi\cdot \widehat{\u}=0,\\
         \mu=-\Delta_\Gamma\vphi+\Psi'(\vphi).
     \end{cases}
     \label{detract}
 \end{align}
Here we have $\V(0)=\V_0:=\v_0-\widehat{\u}(0)$ and $\vphi(0)=\vphi_0$. We will exploit this equivalent problem in the next step.

\textbf{Step 2. Extension of $T_0$ to $T$.} We need to show that actually there exists a strong solution which is globally defined over $[0,T]$, by means of a contradiction argument using energy estimates. Assume then by contradiction that there exists a maximal time of existence $T_m\in[T_0,T)$.  This means that there is no strict extension of $(\v, p,\vphi)$ defined on a strictly larger time interval $[0,T_m')$ being also a local strong solution therein. 

Then it holds, for any $0<\widehat{T}<T_m$, 
\begin{align*}
	&\mathbf{\v} \in  L ^2_{\H^2(\widehat{T})}\cap H^1_{\L^2_\sigma(\widehat{T})} ,\quad p\in L^2_{H^1(\widehat{T})},\quad
	\vphi\in L^\infty_{W^{2,q}(\widehat{T})} \cap H^1_{{H}^1(\widehat{T})} ,\quad \forall q\in[2,+\infty),\\&
\norma{\varphi}<1\quad\text{a.e. in }\gam(t) \quad\text{for a.a. }t\in[0,T_m), \quad
	\mu\in {L}^\infty_{H^1(\widehat{T})}\cap L^2_{H^3(\widehat{T})}.
	\end{align*}
 Now we consider the corresponding triple $(\V,p,\vphi)$, satisfying \eqref{detract}, i.e., we set $\V=\v-\widehat{\u}$. This divergence-free velocity enjoys the same regularity properties as $\v$.
Now we show that all the norms in the spaces above are uniformly bounded by a constant in $(0,T_m)$ and this allows to define the solution also at $\widehat{T}=T_m$. First, let us observe that the regularity above is enough to perform all the following computations, which are therefore rigorous in this framework. We divide the estimates in many steps. All the estimates are related to problem \eqref{detract}.

\textit{First estimate: energy identity. }We primarily notice that it holds, recalling the definition of $\mu$,
\begin{align}
  \nonumber  -&\P\divg(\nablag\vphi\otimes \nabla\vphi)=-\nablag\vphi\Delta_{\gam}\vphi-\nablag(\nablag\vphi)\nablag\vphi\\&=\nablag\vphi(\mu-\Psi'(\vphi))-\frac12\nablag(\vert \nablag\vphi\vert^2)
    =\mu\nablag\vphi-\nablag\left(\Psi(\vphi)-\frac12\vert \nablag\vphi\vert^2\right).
    \label{rewrite}
\end{align}
Note that, in the first identity, we exploited 
\begin{align*}
    \divg(\nablag\vphi\otimes \nabla\vphi)=-\nablag\vphi\Delta\vphi-\widehat{\nabla}_\Gamma(\nablag\vphi)\nablag\vphi,
\end{align*}
and then observed that (see, for instance, \cite[Lemma 15]{BGNhandbook})
\begin{align*}
\widehat{\nabla}_\Gamma(\nablag\vphi)\nablag\vphi=\frac12 \nablag\norma{\nablag\vphi}^2-[\nablag\vphi\cdot((\nablag\mathbf n)\nablag\vphi)]\mathbf n+(\mathbf n\cdot \nablag \vphi)((\nablag\mathbf n)\nablag\vphi),
\end{align*}
so that 
$$
{\nabla}_\Gamma(\nablag\vphi)\nablag\vphi=\P\widehat{\nabla}_\Gamma(\nablag\vphi)\nablag\vphi=\frac12 \nablag\norma{\nablag\vphi}^2.
$$
Therefore, up to redefining the pressure $\tilde{p}=p+\Psi(\vphi)+\frac 12\vert \nablag\vphi\vert^2$ we can substitute the quantity $ -\P\divg(\nablag\vphi\otimes \nabla\vphi)$ with $\mu\nablag\vphi$. 
Again, in the estimates we will consider system \eqref{detract}, and then recover the original velocity. 
Note that, by \eqref{dt1}, for any $t\in[0,T_m]$,
\begin{align}
    \frac12\ddt \intg \rho\norma{\V}^2\ds=\frac12 \intg \dt\rho \norma{\V}^2\ds+\intg\rho \V\cdot \dt\V \ds+\frac 12 \intg \rho\norma{\V}^2Hv_\n\ds.
    \label{rhoas}
\end{align}
Moreover, after integration by parts, recalling that $\J_\rho\cdot \n=0$,
\begin{align}
 0=\intg \divg \left(\frac12\norma{\V}^2\J_\rho\right) \ds=\intg (\J_\rho\cdot\nablag)\V\cdot \V \ds+\intg \frac12 \norma{\V}^2\divg (\J_\rho)\ds, 
 \label{esta1}
\end{align}
as well as, since $\V\cdot \n=0$ and $\divg\V=-\divg\widehat{\u}=-Hv_\n$ by construction,
\begin{align}  &\nonumber0=\intg\divg\left(\rho\V\frac12\norma{\V}^2\right)\ds=\intg \V\cdot \nablag\rho\frac12\norma{\V}^2\ds+\intg \rho \frac12\norma{\V}^2\divg\V\ds+\intg \rho (\V\cdot\nablag)\V\cdot\V\ds\\&
=\intg \V\cdot \nablag\rho\frac12\norma{\V}^2\ds-\intg\frac12\rho\norma{\V}^2Hv_\n\ds+\intg \rho (\V\cdot \nablag)\V\cdot\V\ds.
\label{esta}
\end{align}
Additionally, we get, recalling that the density $\rho$ satisfies (see \eqref{ma}), since $\divg(\V+\widehat{\u}+\VVn)=0$,
$$
0=\dt\rho+\divg(\J_\rho)+\divg( \rho(\V+\widehat{\u}+\VVn))=\dt\rho+\divg(\J_\rho)+ \nablag\rho\cdot (\V+\widehat{\u}),
$$
we deduce, after an integration by parts on the term with $\widehat{ \u}$, noting that $\divg\widehat{\u}=Hv_\n$,
\begin{align*}
&\frac12 \intg \dt\rho \norma{\V}^2\ds\\&=-\intg \divg(\J_\rho)\frac12\norma{\V}^2\ds-\intg \frac12\norma{\V}^2\nablag\rho\cdot \V\ds+\intg \frac12\rho\norma{\V}^2\rho Hv_\n\ds+\intg \rho(\nablag\V)\widehat{\u}\cdot \V \ds.
\end{align*}
Then, from \eqref{esta1} and \eqref{esta}, we deduce
\begin{align*}
    \frac12 \intg \dt\rho \norma{\V}^2\ds= \intg (\nablag\V)(\J_\rho+\rho(\V+\widehat{\u}))\cdot \V\ds.
\end{align*}
Exploiting \eqref{detract}$_1$ and \eqref{rewrite}, we also deduce, integrating by parts,
\begin{align*}
    \intg \rho \V\cdot\dt\V\ds &=-\intg ((\J_\rho+\rho(\V+\widehat{\u}))\cdot\nablag)\V)\cdot \V\ds-\intg \rho (\V\cdot \nablag)\widehat{\u}\cdot \V\ds\\&-\intg \rho v_\n \H{\V}\cdot\V\ds-\intg v_\n\H\J_\rho\cdot \V-2\intg \nu(\vphi)\norma{\E_S(\V)}^2\ds\\&+\intg \mu\nablag\vphi\cdot \V\ds-2\intg  \nu(\vphi)v_\n\H :\nablag\V\ds+\frac12\intg \rho\nablag (v_\n)^2\cdot \V\ds\\& 
    -\intg \rho \dt\widehat{\u}\cdot \V-((\rho\widehat{\u}+\J_\rho)\cdot \nablag)\widehat{\u}\cdot \V\ds\\&-\intg \rho v_\n\H\widehat{\u}\cdot \V\ds-2\intg \nu(\vphi)\E_S(\widehat{\u}):\E_S(\V)\ds.
\end{align*}
Putting all these identities above in \eqref{rhoas}, we infer
\begin{align*}
        &\frac12\ddt\intg \rho\norma{\V}^2\ds+2\intg \nu(\vphi)\norma{\E_S(\V)}^2\ds\\&=\frac12\intg\rho\norma{\V}^2H v_\n\ds-\intg \rho (\V\cdot \nablag)\widehat{\u}\cdot \V\ds-\intg \rho v_\n \H{\V}\cdot\V\ds-\intg v_\n\H\J_\rho\cdot \V\ds\\&
        +\intg \mu\nablag\vphi\cdot \V\ds-2\intg  \nu(\vphi)v_\n\H : \nablag\V\ds+\frac12\intg \rho\nablag (v_\n)^2\cdot \V\ds\\& 
    -\intg \rho \dt\widehat{\u}\cdot \V-((\rho\widehat{\u}-\J_\rho)\cdot \nablag)\widehat{\u}\cdot \V\ds-\intg \rho v_\n\H\widehat{\u}\cdot \V\ds-2\intg \nu(\vphi)\E_S(\widehat{\u}):\E_S(\V)\ds.
\end{align*}
From \eqref{detract}$_3$, multiplying by $\mu$, integrating by parts where necessary, we obtain (see also \eqref{dtfin} with $\v=\V+\widehat{\u}+\VVn$)
\begin{align}
        \nonumber&\frac{d}{dt}\left(\frac12\Vert \nablag\vphi\Vert^2+\intg\Psi(\vphi)\ds\right)+\intg\vert \nablag\mu\vert^2\ds+\intg \nablag\vphi\cdot (\widehat{\u}+\V) \ \mu\ds\\&=\frac12\intg \vert \nablag\vphi\vert^2H v_\n\ds -\intg \nablag\vphi \cdot\E_S(\V_\n)\nablag\vphi\ds+\intg\Psi(\vphi)H v_\n\ds,
        \label{dtfin2}
    \end{align}
since $\nablag\vphi\cdot \VVn=0$. Summing up the two identities above, and setting
 \begin{align}
E_{tot}(\rho,\V,\vphi):=\frac 12\intg \rho\norma{\V}^2\ds+ \frac12\Vert \nablag\vphi\Vert^2+\intg\Psi(\vphi)\ds,
    \label{Energ}
\end{align}
we get the following energy identity
\begin{align}
&\nonumber\ddt E_{tot}(\rho,\V,\vphi)+2\intg \nu(\vphi)\norma{\E_S(\V)}^2\ds+\intg\vert \nablag\mu\vert^2\ds\\&\nonumber=\frac12\intg\rho\norma{\V}^2H v_\n\ds-\intg \rho (\V\cdot \nablag)\widehat{\u}\cdot \V\ds-\intg \rho v_\n \H{\V}\cdot\V\ds-\intg v_\n\H\J_\rho\cdot \V\ds\\&\nonumber
        -\intg \mu\nablag\vphi\cdot \widehat{\u}\ds-2\intg  \nu(\vphi)v_\n\H : \nablag\V\ds+\frac12\intg \rho\nablag (v_\n)^2\cdot \V\ds\\& 
    -\intg \rho \dt\widehat{\u}\cdot \V-((\rho\widehat{\u}-\J_\rho)\cdot \nablag)\widehat{\u}\cdot \V\ds-\intg \rho v_\n\H\widehat{\u}\cdot \V\ds-2\intg \nu(\vphi)\E_S(\widehat{\u}):\E_S(\V)\ds\nonumber\\&
    +\frac12\intg \vert \nablag\vphi\vert^2H v_\n\ds -\intg \nablag\vphi\cdot \E_S(\V_\n)\nablag\vphi\ds+\intg\Psi(\vphi)H v_\n\ds,
    \label{energyfinal}
\end{align}
since, due to $\divg\V=0$, integrating by parts, $\frac12\intg \nablag (v_\n)^2\cdot \V\ds=0$.
Now we recall that, thanks to the regularity of the flow map (see \eqref{regularity_basic}-\eqref{regularity_basic2}), we have
\begin{align}
&\nonumber\norm{\H}_{\mathbf C^4\gt}+\norm{\VVn}_{\mathbf C^4\gt}+\norm{\A}_{C^4\gt}\\&+\norm{\A^{-1}}_{C^4\gt}+\norm{\dt\A^{-1}}_{C^3\gt}+\norm{\widehat{\u}}_{C^4_{\W^{3,p}}}\leq C(T),\quad \forall t\in[0,T],\quad \forall p\in[2,+\infty),  \label{regflowmap}
\end{align}
and we recall that, for $p>2$, $\W^{3,p}\gt\hookrightarrow \W^{2,\infty}\gt$.
We can thus estimate all the terms accordingly. In particular, we have, recalling the definition of $\J_\rho$, and by H\"{o}lder's and Young's inequality, 
\begin{align*}
    &\norma{\frac12\intg\rho\norma{\V}^2H v_\n\ds-\intg \rho (\V\cdot \nablag)\widehat{\u}\cdot \V\ds-\intg \rho v_\n \H{\V}\cdot\V\ds-\intg v_\n\H\J_\rho\cdot \V\ds}\\&
    \leq C\norm{Hv_\n}_{L^\infty\gt}\norm{\V}^2+C\norm{\V}^2\norm{\nablag\widehat{\u}}_{\L^\infty\gt}+C\norm{v_\n\H}_{\L^\infty\gt}\norm{\V}^2+C\norm{v_\n\H}_{\L^\infty\gt}\norm{\nablag\mu}\norm{\V}\\&
    \leq C(T)\norm{\V}^2+\frac14 \norm{\nablag\mu}^2.
\end{align*}
Then, after an integration by parts, recalling $\divg\widehat{\u}=Hv_\n$, we infer
\begin{align*}
    \norma{-\intg \mu\nablag\vphi\cdot \widehat{\u}\ds}\leq \norma{\intg \vphi\nablag\mu\cdot \widehat{\u}\ds}+ \norma{\intg \mu\vphi Hv_\n\ds}.
\end{align*}
By following the same argument leading to \eqref{mubar}, recalling that $\norma{\vphi}<1$ on $[0,T_m)$ by assumption, we deduce that
\begin{align}
\norma{\overline{\mu}}\leq C(T)(1+\norm{\nablag\mu})
    \label{mub}
\end{align}
so that, since $\norma{\Gamma(t)}\equiv \norma{\gz}$, 
\begin{align*}
    \norma{-\intg \mu\nablag\vphi\cdot \widehat{\u}\ds}&\leq C\norm{\vphi}_{L^\infty\gt}\norm{\mu}\norm{\widehat{\u}}_{\L^2\gt}
    \leq C(T)(1+\norm{\nablag\mu})\leq C(T)+\frac14 \norm{\nablag\mu}^2.
\end{align*}
Then, we have, recalling again the definition of $\J_\rho$ and Young's inequality,
\begin{align*}
&    \left\vert-2\intg  \nu(\vphi)v_\n\H:\nablag\V\ds
    -\intg \rho \dt\widehat{\u}\cdot \V-((\rho\widehat{\u}-\J_\rho)\cdot \nablag)\widehat{\u}\cdot \V\ds+\frac12\intg \rho\nablag (v_\n)^2\cdot \V\ds\right.\\&\left.-\intg \rho v_\n\H\widehat{\u}\cdot \V\ds-2\intg \nu(\vphi)\E_S(\widehat{\u}):\E_S(\V)\ds\right\vert
    \\&
    \leq C\norm{v_\n\H}_{\L^\infty\gt}\norm{\nablag\V}+C\norm{\dt\widehat{\u}}_{\L^\infty\gt}\norm{\V}+C\norm{\widehat{\u}}_{\L^\infty\gt}\norm{\nablag\widehat{\u}}_{\L^\infty\gt}\norm{\V}\\&
    +C\norm{\nablag\mu}\norm{\nablag\widehat{\u}}_{\L^\infty\gt}\norm{\V}+C\norm{\nablag(\vn)^2}_{\Linfg}\norm{\V}\\&+C\norm{v_\n\H\widehat{\u}}_{\L^\infty\gt}\norm{\V}+C\norm{\widehat{\u}}_{\W^{1,\infty}\gt}\norm{\nablag\V}\\&
    \leq C(T)(1+\norm{\V}^2)+\frac14 \norm{\nablag\mu}^2+{\nu_*}\norm{\E_S(\V)}^2,
\end{align*}
where in the last estimate we applied the (evolving) Korn's inequality \eqref{Korn}.
By similar estimates, we infer 
\begin{align*}
    &\norma{\frac12\intg \vert \nablag\vphi\vert^2H v_\n\ds -\intg \nablag\vphi \E_S(\V_\n)\nablag^T\vphi\ds+\intg\Psi(\vphi)H v_\n\ds}\\&\leq C(T)\left(\norm{\nablag\vphi}^2+\intg (\Psi(\vphi)+\widehat{C})\ds\right),
\end{align*}
where $\widehat{C}>0$ is such that $\Psi+\widehat{C}>0$. By redefining $E_{tot}$ as to be 
\begin{align}
\widehat{E}_{tot}(t):=\frac 12\intg \rho\norma{\V}^2\ds+ \frac12\Vert \nablag\vphi\Vert^2+\intg(\Psi(\vphi)+\widehat{C})\ds
    \label{Energ2},
\end{align}
and recalling that $\nu(\vphi)\geq \nu_*$ by assumption,
we deduce from all the estimates above put together in \eqref{energyfinal},
that
\begin{align}
    \ddt \widehat{E}_{tot}(t)+\nu_*\norm{\E_S(\V)}^2+\frac{1}{4}\norm{\nablag\mu}^2\leq C(T)(1+\widehat{E}_{tot}(t)), \quad \forall t\in[0,T_m).
    \label{etot}
\end{align}
Since $T_m\leq T$, this allows to deduce by the Gronwall lemma and \eqref{mub} (by Poincaré's inequality), that
\begin{align}
    \norm{\V}_ {L^\infty_{\L^2_\sigma(T_m)}}+\norm{\V}_{L^2_{\H^1_\sigma(T_m)}}+\norm{\vphi}_{L^\infty_{H^1(T_m)}}+\norm{\mu}_{L^2_{H^1(T_m)}}\leq C(T).
    \label{uniformtime}
\end{align}
Since $\vphi$ satisfies problem \eqref{pr} with $\v=\V+\widehat{\u}+\VVn$, the regularity on the initial datum $\vphi_0$, together with estimates \eqref{regflowmap} and \eqref{uniformtime} allow to apply Theorem \ref{existence} and deduce, by uniqueness of the solution to \eqref{pr}, that  
 \begin{align} \norm{\mu}_{H^1_{H^{-1}(T_m)}}+\norm{\mu}_{L^\infty_{H^1(T_m)}}+\norm{\mu}_{L^2_{H^3(T_m)}}+\norm{\vphi}_{L^\infty_{H^3(T_m)}}+\norm{\vphi}_{H^1_{H^1(T_m)}}\leq C(T,\norm{\u}_{L^2_{\H^1(T_m)
 }})\leq C(T),
 \label{newreg2}
 \end{align}
as well as there exists $\delta=\delta(T)\in(0,1)$ such that 
\begin{align}
\sup_{t\in[0,T_m)}\norm{\vphi}_{L^\infty\gt}\leq 1-\delta.
\label{separation_ess}
\end{align}
\textit{Higher-order estimates. }To conclude the contradiction argument, we need to show some higher-order estimates on the velocity $\V$. First, we notice that \eqref{detract}$_1$, satisfied by $(\V,p,\vphi)$, can be rewritten as 

\begin{align*}
   &-2\P\divg(\nu(\vphi)\E_S(\V))+\omega\V+\nablag \pi
     \\&=\mathbf f:= \omega\V-\rho \P\dt \V -((\rho\V+\J_\rho)\cdot\nabla_\Gamma)\V-\rho(\widehat{\u}\cdot\nabla_\Gamma)\V-\rho(\V\cdot\nabla_\Gamma)\widehat{\u}-\rho \Vn\mathbf{H}\V-\Vn\H\J_\rho\\&\quad -\P\divg(\nablag\vphi\otimes\nablag\vphi)+2\P\divg(\nu(\vphi)\Vn\mathbf{H})
     +\frac\rho2\nablag^T(\Vn)^2-\rho \P\dt \widehat{\u} -((\rho\widehat{\u}+\J_\rho)\cdot\nabla_\Gamma)\widehat{\u}\\&\quad-\rho \Vn\mathbf{H}\widehat{\u}+2\P\divg(\nu(\vphi)\E_S(\widehat{\u})),
\end{align*}
where $\omega>0$ is arbitrary. This means that the corresponding weak formulation coincides with the one in \cite[Lemma 7.4]{AGP1}. Therefore, we can apply \cite[Lemma 7.4]{AGP1} and infer that
\begin{align}
\label{VVV}
\norm{\V(t)}_{{\H^2}(t)}\leq C(T,\omega)(1+\norm{\vphi(t)}_{L_{W^{1,\infty}\gt}})\norm{\f(t)}_{{\L^2}\gt}\leq C(T,\omega)\norm{\f(t)}_{{\L^2}\gt},
\end{align}
for almost any $t\in(0,T_m)$, recalling the embedding (with a uniform-in-time embedding constant) $H^3\gt\hookrightarrow W^{1,\infty}\gt$ and \eqref{newreg2}. Now we need to estimate the right-hand side. We have, recalling \eqref{uniformtime},
$$
\norm{\omega\V-\rho \P\dt \V}\leq C\norm{\V}+\rho^*\norm{\dt\V})\leq C(T)+\rho^*\norm{\P\dt\V},
$$
as well as, by Sobolev-Gagliardo-Nirenberg's and Young's inequalities, and \eqref{Korn},
\begin{align}
    \nonumber\norm{\rho(\V\cdot\nablag)\V}&\leq \rho^*\norm{\V}_{\L^4\gt}\norm{\nablag\V}_{\L^4\gt}\\&\nonumber\leq C(T)\norm{\V}^\frac12\norm{\V}_{\H^1\gt}\norm{\V}_{\H^2\gt}^\frac12\\&\nonumber
    \leq C(T)\norm{\V}^\frac12(\norm{\V}+\norm{\E_S(\V)})\norm{\V}_{\H^2\gt}^\frac12\\&
    \nonumber\leq 
    C(T)(1+\norm{\E_S(\V)})\norm{\V}_{\H^2\gt}^\frac12\\&
    \leq C(T)(1+\norm{\E_S(\V)}^2)+\frac{1}{2C(T,\omega)}\norm{\V}_{\H^2\gt}.\label{convective}
\end{align}
Then, by Agmon's inequality \eqref{Agmon} and \eqref{uniformtime}-\eqref{newreg2}, we infer, recalling the definition of $\J_\rho$, 
\begin{align}
    \nonumber\norm{(\J_\rho\cdot\nablag)\V}&\leq C\norm{\nablag\mu}_{\L^\infty\gt}\norm{\nablag\V}\\&\nonumber\leq C(T)\norm{\nablag\mu}^\frac12\norm{\mu}_{H^3\gt}^\frac12(\norm{\V}+\norm{\E_S(\V)})\\&
    \leq C(T)\norm{\mu}_{H^3\gt}^\frac12(1+\norm{\E_S(\V)}).
    \label{Jrro}
\end{align}
Furthermore, again by \eqref{newreg2}, we deduce 
\begin{align*}
\norm{\P\divg(\nablag\vphi\otimes\nablag\vphi)}&\leq C\norm{\vphi}_{H^2\gt}\norm{\nablag\vphi}_{\L^\infty\gt}
\leq C(T),
\end{align*}
again by the embedding $H^3\gt\hookrightarrow W^{1,\infty}\gt$.
Analogously, we have
\begin{align*}
&\norm{2\P\divg(\nu(\vphi)\Vn\mathbf{H})+2\P\divg(\nu(\vphi)\E_S(\widehat{\u}))}\\&\leq C(1+\norm{\vphi}_{H^1\gt})(\norm{v_\n\H}_{\W^{1,\infty}\gt}+\norm{\widehat{\u}}_{\W^{2,\infty}\gt})\leq C(T).
\end{align*}
Proceeding in estimating the lower-order terms in $\f$, we have, recalling \eqref{regflowmap}, \eqref{uniformtime}-\eqref{newreg2} and the definition of $\J_\rho$,
\begin{align*}
   &\left\Vert -\rho(\V\cdot\nabla_\Gamma)\widehat{\u}-\rho \Vn\mathbf{H}\V-\Vn\H\J_\rho
     +\frac\rho2\nablag^T(\Vn)^2-\rho \P\dt \widehat{\u} -((\rho\widehat{\u}+\J_\rho)\cdot\nabla_\Gamma)\widehat{\u}-\rho \Vn\mathbf{H}\widehat{\u}\right\Vert\\&
     \leq C(T)(\norm{\V}\norm{\nablag\widehat{\u}}_{\Linfg}+\norm{v_\n\H}_{\Linfg}(\norm{\V}+\norm{\nablag\mu})+\norm{\nablag(v_\n)^2}_{\Linfg}\\&+\norm{\dt\widehat{\u}}_{\Linfg}+(\norm{\widehat{\u}}_{\Linfg}+\norm{\nablag\mu})\norm{\nablag\wu}_{\Linfg}+\norm{v_\n\H}_{\Linfg})\\&
     \leq C(T).
\end{align*}
Therefore, we can conclude from \eqref{VVV}, elevating to the square, that
\begin{align}
    &\norm{\V}_{\H^2\gt}^2\leq C(T,\omega)(1+\norm{\dt\V}^2+\norm{\E_S(\V)}^4+\norm{\mu}_{H^3\gt}(1+\norm{\E_S(\V)}^2)).
    \label{final1}
\end{align}
We recall that, thanks to \eqref{regflowmap} and \eqref{uniformtime},
\begin{align}
    \norm{\dt\V}\leq C(T)(1+\norm{\dts\V}),\quad \norm{\dts\V}\leq C(T)(1+\norm{\dt\V}),\label{equivalence1}
\end{align}
so that 
\begin{align}
    &\norm{\V}_{\H^2\gt}^2\leq C(T,\omega)(1+\norm{\dts\V}^2+\norm{\E_S(\V)}^4+\norm{\mu}_{H^3\gt}(1+\norm{\E_S(\V)}^2)).
    \label{final2}
\end{align}
Let us now multiply \eqref{detract}$_1$ by $\dts\V$ and integrate over $\Gamma(t)$. Recall that $\divg(\dts\V)=0$ and $\P\dts\V=\dts\V$. We cannot simply use $\dt\V$, since it does not preserve the divergence-free property. The most involved part is related to the following term 
$$-\intg \divg(\nu(\vphi)\E_S(\V))\cdot \dts\V\ds.$$
Indeed, from \eqref{Ess} and Remark \ref{basic_rem} we infer 
\begin{align}
&\nonumber\frac{d}{dt}\intg \nu(\vphi) \norma{\E_S(\V)}^2\ds=\int_{\gam}\nu'(\vphi)\dt\vphi\norma{\E_S(\V)}^2\ds\nonumber-2\int_\gam \divg(\nu(\vphi)\E_S( \V))\cdot\dt\V\ds\\&+4\int_{\gam}\nu(\vphi)\mathcal{S}((\n\otimes \n)\widehat{\nabla}_\gam\V_\n)\widehat{\E}_S(\V):\E_S(\V)\ds\nonumber+4\int_{\gam}\nu(\vphi){\E}_S(\V):\widehat{\E}_S(\V)\mathcal{S}((\n\otimes \n)\widehat{\nabla}_\gam\V_\n)\ds\non\\&-2\int_\gam \nu(\vphi)\mathcal{S}((\nablag\V)\nablag\V_\n):\E_S(\V)\ds+\intg \nu(\vphi) \norma{\E_S(\V)}^2\ \divg\V_\n\ds.
\label{Ess1}
\end{align}
The difficulty is that we only control $\dts\V$ (or $\P\dt\V$) and not the entire $\dt\V$. 
Recalling \eqref{relation_base1}, we perform the following splitting
\begin{align*}
    &-2\int_\gam \divg(\nu(\vphi)\E_S( \V))\cdot\dts\V\ds\\&=-2\int_\gam \divg(\nu(\vphi)\E_S( \V))\cdot\dt\V\ds+2\int_\gam \divg(\nu(\vphi)\E_S( \V))\cdot(\A(\dt\A^{-1})\V)\ds.
\end{align*}

In the end we have:
\begin{align}
  &\nonumber \frac{d}{dt}\intg \nu(\vphi) \norma{\E_S(\V)}^2\ds+\intg\rho\norma{\dts\V}^2\ds\\&=-\nonumber\intg\rho\dts\V\cdot \A(\dt\A^{-1})\V\ds+\int_{\gam}\nu'(\vphi)\dt\vphi\norma{\E_S(\V)}^2\ds-2\int_\gam \divg(\nu(\vphi)\E_S( \V))\cdot\A(\dt\A^{-1})\V\ds\\&\nonumber+4\int_{\gam}\nu(\vphi)\mathcal{S}((\n\otimes \n)\widehat{\nabla}_\gam\V_\n)\widehat{\E}_S(\V):\E_S(\V)\ds\nonumber+4\int_{\gam}\nu(\vphi){\E}_S(\V):\widehat{\E}_S(\V)\mathcal{S}((\n\otimes \n)\widehat{\nabla}_\gam\V_\n)\ds\\&\nonumber-2\int_\gam \nu(\vphi)\mathcal{S}((\nablag\V)\nablag\V_\n):\E_S(\V)\ds\nonumber+\intg \nu(\vphi) \norma{\E_S(\V)}^2\ \divg\V_\n\ds\\&\nonumber
   -\intg ((\rho\V+\J_\rho)\cdot\nabla_\Gamma)\V\cdot \dts\V\ds-\intg \rho(\widehat{\u}\cdot\nabla_\Gamma)\V\cdot \dts\V\ds\\&\nonumber+\intg\left(-\rho(\V\cdot\nabla_\Gamma)\widehat{\u}-\rho \Vn\mathbf{H}\V-\Vn\H\J_\rho\right)\cdot \dts\V\ds\\&\nonumber+\intg\left(-\P\divg(\nablag\vphi\otimes\nablag\vphi)+2\P\divg(\nu(\vphi)\Vn\mathbf{H})\right)\cdot \dts\V\ds\\&
     +\intg\left(\frac\rho2\nablag(\Vn)^2-\rho \P\dt \widehat{\u} -((\rho\widehat{\u}+\J_\rho)\cdot\nabla_\Gamma)\widehat{\u}-\rho \Vn\mathbf{H}\widehat{\u}+2\P\divg(\nu(\vphi)\E_S(\widehat{\u}))\right)\cdot \dts\V\ds.
     \label{strongest}
\end{align}
We now estimate all the terms appearing here.
First, we have from \eqref{regflowmap} (see also \eqref{regularity_basic}) and \eqref{uniformtime},
\begin{align*}
    \norma{-\intg\rho\dts\V\cdot \A(\dt\A^{-1})\V\ds}\leq C(T)\norm{\dts\V}\norm{\V}\leq C(T)+\frac{\rho_*}{16}\norm{\dts\V}^2.
\end{align*}
Then, recalling \eqref{Korn}, \eqref{Gagliardo}, and \eqref{newreg2},
\begin{align*} &\norma{\int_{\gam}\nu'(\vphi)\dt\vphi\norma{\E_S(\V)}^2\ds}\leq C\norm{\dt\vphi}\norm{\E_S(\V)}^2_{\L^4(\gt)}\\&
\leq C\norm{\dt\vphi}\norm{\E_S(\V)}\norm{\V}_{\H^2(\Gamma(t))}
\leq C(T,\epsilon)\norm{\dt\vphi}^2\norm{\E_S(\V)}^2+\epsilon\norm{\V}_{\H^2\gt}^2,
\end{align*}
for some $\eps>0$ suitably chosen later on. Then, thanks to \eqref{regflowmap} and \eqref{newreg2},
\begin{align*}
    &\norma{-2\int_\gam \divg(\nu(\vphi)\E_S( \V))\cdot \A(\dt\A^{-1})\V\ds}\leq C\norm{\vphi}_{W^{1,\infty}\gt}\norm{\E_S( \V)}\norm{\V}+C\norm{\V}_{\H^2\gt}\norm{\V}\\&
    \leq C(T)\norm{\E_S(\V)}+C(T)\norm{\V}_{\H^2\gt}\leq
    C(T,\eps)(1+\norm{\E_S(\V)}^2)+\eps\norm{\V}_{\H^2\gt}^2.
\end{align*}
Proceeding in the estimates, recalling \eqref{Korn} and \eqref{regflowmap}, we have
\begin{align*}
    &\left\vert 4\int_{\gam}\nu(\vphi)\mathcal{S}((\n\otimes \n)\widehat{\nabla}_\gam\V_\n)\widehat{\E}_S(\V):\E_S(\V)\ds\nonumber+4\int_{\gam}\nu(\vphi){\E}_S(\V):\widehat{\E}_S(\V)\mathcal{S}((\n\otimes \n)\widehat{\nabla}_\gam\V_\n)\ds\right.\\&\left.-2\int_\gam \nu(\vphi)\mathcal{S}((\nablag\V)\nablag\V_\n):\E_S(\V)\ds\nonumber+\intg \nu(\vphi) \norma{\E_S(\V)}^2\ \divg\V_\n\ds\right\vert\\&
    \leq C(T)\norm{\V}_{\H^1\gt}^2\leq C(T)(\norm{\V}^2+\norm{\E_S(\V)}^2)\leq C(T)(1+\norm{\E_S(\V)}^2).
\end{align*}
Moreover, by similar estimates as in \eqref{convective}, we get
\begin{align*}
\norma{\intg\rho(\V\cdot\nablag)\V\cdot \dts\V\ds}&\leq \rho^*\nonumber\norm{(\V\cdot\nablag)\V}\norm{\dts\V}\leq \rho^*\norm{\V}_{\L^4\gt}\norm{\nablag\V}_{\L^4\gt}\norm{\dts\V}\\&
    \nonumber\leq 
    C(T)(1+\norm{\E_S(\V)})\norm{\V}_{\H^2\gt}^\frac12\norm{\dts\V}\\&
    \leq C(T,\eps)(1+\norm{\E_S(\V)}^4)+\frac{\rho_*}{16}\norm{\dts\V}^2+\eps\norm{\V}_{\H^2\gt}^2,
\end{align*}
and, similarly to \eqref{Jrro}, \begin{align*}
&\norma{\intg (\J_\rho\cdot\nablag)\V\cdot \dts\V\ds}\leq \nonumber\norm{(\J_\rho\cdot\nablag)\V}\norm{\dts\V}\leq C\norm{\nablag\mu}_{\L^\infty\gt}\norm{\nablag\V}\norm{\dts\V}\\&\nonumber\leq C(T)\norm{\nablag\mu}^\frac12\norm{\mu}_{H^3\gt}^\frac12(\norm{\V}+\norm{\E_S(\V)})\norm{\dts\V}
    \leq C(T)\norm{\mu}_{H^3\gt}(1+\norm{\E_S(\V)}^2)+\frac{\rho_*}{16}\norm{\dts\V}^2.
\end{align*}
We can then proceed in the estimates: recalling that $\divg\dts\v=0$ and $\P\dts\v=\dts\v$, as well as \eqref{rewrite} and \eqref{newreg2}, we get   
\begin{align*}
&\norma{\intg \P\divg(\nablag\vphi\otimes\nablag\vphi)\cdot \dts\V\ds}=\norma{\intg \mu\nablag\vphi\cdot \dts\V\ds}
\leq \norma{\intg \nablag\mu\vphi\cdot \dts\V\ds}\\&
\leq \norm{\nablag\mu}\norm{\vphi}_{L^\infty\gt}\norm{\dts\V}
\leq C(T)+\frac{\rho_*}{16}\norm{\dts\V}^2.
\end{align*}
Concerning the lower order terms, we then observe that, thanks to \eqref{Korn}, and \eqref{uniformtime}, \eqref{newreg2}, and recalling the definition of $\J_\rho$,
\begin{align*}
  & \left\vert -\intg \rho(\widehat{\u}\cdot\nabla_\Gamma)\V\cdot \dts\V\ds+\intg\left(-\rho(\V\cdot\nabla_\Gamma)\widehat{\u}-\rho \Vn\mathbf{H}\V-\Vn\H\J_\rho\right)\cdot \dts\V\ds\right\vert\\&
  \leq C(\norm{\widehat{\u}}_{\L^\infty\gt}\norm{\nablag\V}\norm{\dts\V}+\norm{\V}\norm{\nablag\wu}_{\Linfg}\norm{\dts\V}+(\norm{\V}+\norm{\nablag\mu})\norm{v_\n\H}_{\Linfg}\norm{\dts\V})\\&\leq C(T)(1+\norm{\E_s(\V)}^2)+\frac{\rho_*}{16}\norm{\dts\V}^2, 
\end{align*}
as well as, 
\begin{align*}
    &\norma{\intg\left(\frac\rho2\nablag(\Vn)^2-\rho \P\dt \widehat{\u} -((\rho\widehat{\u}+\J_\rho)\cdot\nabla_\Gamma)\widehat{\u}-\rho \Vn\mathbf{H}\widehat{\u}\right)\cdot \dts\V\ds}\\&
    \leq C(T)\left(\norm{\nablag (v_\n)^2}_{\Linfg}+\norm{\dt\wu}_{\Linfg}\right.\\&\left.+(\norm{\wu}_{\Linfg}+\norm{\nablag\mu})\norm{\nablag\wu}_{\Linfg}+C\norm{v_\n\H}_{\Linfg}\norm{\wu}_{\Linfg}\right)\norm{\dts\V}\\&
    \leq C(T)+\frac{\rho^*}{16}\norm{\dts\V}^2.
\end{align*}
In conclusion, we have
\begin{align*}
&\norma{\intg\left(2\P\divg(\nu(\vphi)\E_S(\widehat{\u}))+2\P\divg(\nu(\vphi)\Vn\mathbf{H})\right)\cdot \dts\V\ds}\\&
\leq C(\norm{\vphi}_{H^1\gt}\norm{\E_S(\wu)}_{\Linfg}+\norm{\wu}_{\H^2\gt})\norm{\dts\V}\\&
+C(\norm{\vphi}_{H^1\gt}\norm{v_\n\H}_{\Linfg}+\norm{v_\n\H}_{\H^1\gt})\norm{\dts\V}\\&
\leq C(T)+\frac{\rho_*}{16}\norm{\dts\V}^2.
\end{align*}
Therefore, to sum up, we have from \eqref{strongest}, recalling that $\rho\geq \rho_*$,
\begin{align*}
    &\frac{d}{dt}\intg \nu(\vphi) \norma{\E_S(\V)}^2\ds+\frac{9\rho_*}{16}\norm{\dts\V}^2\\&\leq C(T,\epsilon)(1+\norm{\dt\vphi}^2+\norm{\mu}_{H^3\gt}+\norm{\E_S(\V)}^2)\norm{\E_S(\V)}^2+3\eps\norm{\V}_{\H^2\gt}^2.
\end{align*}
In conclusion, we sum up this inequality with \eqref{final1} multiplied by $\tfrac{\rho_*}{16C(\omega,T)}$, and we choose $\eps=\tfrac{\rho_*}{64C(\omega,T)}$, to deduce in the end that
\begin{align*}
   & \frac{d}{dt}\intg \nu(\vphi) \norma{\E_S(\V)}^2\ds+\frac{\rho_*}{2}\norm{\dts\V}^2+\tfrac{\rho_*}{64C(\omega,T)}\norm{\V}_{\H^2\gt}^2\\&\leq C(T)(1+\norm{\dt\vphi}^2+\norm{\mu}_{H^3\gt}+\norm{\E_S(\V)}^2)\norm{\E_S(\V)}^2,\quad \text{for almost any }t\in(0,T_m),
\end{align*}
which entails, thanks to \eqref{uniformtime}-\eqref{newreg2} and \eqref{relation_base1}, that, by the Gronwall lemma, since $T_m<T$,
\begin{align}
\norm{\V}_{L^\infty_{\H^1_\sigma(T_m)}}+\norm{\V}_{L^2_{\H^2(T_m)}}+\norm{\V}_{H^1_{\L^2(T_m)}}\leq C(T).
    \label{regv}
\end{align}
Now, let us observe that the regularity \eqref{newreg2}, together with \eqref{regv}, entails (see also Lemma \ref{continuous_embeddings}) that 
$$
\norm{\V}_{C^0_{\H^1_\sigma(T_m)}}+\norm{\vphi}_{C^0_{H^2
(T_m)}}+\norm{\mu}_{C^0_{H^1(T_m)}}\leq C(T).
$$
This means that we can also extend the solution at $t=T_m$, and, in particular, it holds 
\begin{align*}
    \V(T_m)\in \H^1_\sigma(\gam(T_m)),\quad \vphi(T_m)\in H^2(\gam(T_m)),\quad \mu(T_m)\in H^1(\gam(T_m)),
\end{align*}
together with $\norm{\vphi(T_m)}_{L^\infty(\gam(T_m))}\leq 1-2\delta_m$ (see \eqref{separation_ess}), for some $\delta_m\in(0,1)$. Moreover, it holds
$$
-\Delta_\gam\vphi(T_m)+\Psi'(\vphi(T_m))=\mu(T_m),\quad \text{for almost any }x\in \Gamma(T_m),
$$
entailing by elliptic regularity, exploiting the separation property, that 
$$
\vphi(T_m)\in H^3(\gam(T_m))\hookrightarrow B_{p,q}^{4-\frac4q}(\gam(T_m)),
$$
for some $q\in(2,4]$ and $p>4$, see also \eqref{H3b}. This means that the pair $(\v(T_m),\vphi(T_m))$, with $\v(T_m):=\V(T_m)+\widehat{\u}(T_m)$, is sufficiently regular to be considered as initial datum in the application of \cite[Theorem 4.1]{AGP1} (together with \cite[Remarks 4.3-4.4] {AGP1}) to deduce, following word by word the argument in Step 0, that there exists a strong solution according to Definition \ref{strong}, defined on an interval $[T_m,T_1]$, with $T_1>T_m$, extending the one on $[0,T_m)$. This contradicts the maximality of  the considered solution on $(0,T_m)$, entailing that $T_m=T$, i.e., any strong solution is globally defined on $[0,T]$, concluding the first part of the proof. We now need to prove the uniqueness of global strong solutions. 

\textbf{Step 3. Continuous dependence and Uniqueness of strong solutions.}
Let us consider two sets of initial data $(\V_{0,i},\vphi_{0,i})$, $i=1,2$, satisfying the assumptions of Theorem \ref{strong1} and $\overline{\vphi}_1=\overline{\vphi}_2$, so that there exist two global strong solutions $(\V_i,\vphi_i,p_i)$ departing from these data. We also define $\rho_i:=\rho(\vphi_i)$. Let us set $\V:=\V_1-\V_2$, $\vphi:=\vphi_1-\vphi_2$, $\pi=\pi_1-\pi_2$ and consider the equations satisfied by these variables. In particular, we multiply the equation for the difference of velocities by $\V$ and integrate over $\gam(t)$, whereas we take the gradient of the Cahn-Hilliard equation (written for the difference $\vphi$), multiply it by $\nablag\vphi$ and again integrate over $\gam(t)$. We have
\begin{align}
\begin{cases}
    \intg \rho_1\dt\V\cdot\V\ds+\intg (\rho_1-\rho_2)\dt\V_2\cd\V\ds+\ig\rho_1(\V\cd\nablag)\V_1\cd\V\ds\\+\ig(\rho_1-\rho_2)(\V_2\cd\nablag)\V_1\cdot \V\ds
    \ig \ro_2(\V_2\cd\nablag)\V\cd\V\ds+\ig(\J_{\ro_1}\cd\nablag)\V\cd\V\ds\\+\ig((\J_{\ro_1}-\J_{\ro_2})\cd\nablag\V_2)\cd \V\ds
+\ig\ro_1(\wu\cd\nablag)\V\cd \V\ds+\ig \dr(\wu\cd\nablag)\V_2\cd\V\ds\\+\ig\rho_1(\V\cd\nablag)\wu\cd\V\ds+\ig\dr(\V_2\cd\nablag)\wu\cd\V\ds+\ig\ro_1\vn\H\V\cd\V\ds\\+\ig\dr\vn\H\V_2\cd\V\ds+\ig\vn\H(\J_{\ro_1}-\J_{\ro_2})\cd\V\ds\\+2\ig\nu(\vphi_1)\norma{\E_S(\V)}^2\ds+2\ig(\nu(\vphi_1)-\nu(\vphi_2))\E_S(\V_2):\E_S(\V)\ds
\\=\ig (\nablag\vphi_1\otimes \nablag\vphi):\nablag\V\ds+\ig (\nablag\vphi\otimes \nablag\vphi_2):\nablag\V\ds\\-2\ig (\nu(\vphi_1)-\nu(\vphi_2))\vn\H:\nablag\V\ds
-\ig \dr\dt\wu\cd\V\ds-\ig\dr(\wu\cd\nablag)\wu\cd\V\ds\\-\ig((\J_{\ro_1}-\J_{\ro_2})\cd\nablag)\wu\cd\V\ds-\ig\dr(\vn\H\wu-\frac12\nablag(\vn)^2)\cd\V\ds\\-2\ig(\nu(\vphi_1)-\nu(\vphi_2))\E_S(\wu):\E_S(\V)\ds,
\\
\\
\ig\nablag\dt\vphi\cd \nablag\vphi\ds-\intg \nablag\Delta(\mu_1-\mu_2)\cd\nablag\vphi\ds\\+\ig\nablag(\nablag(\vphi_1)\cd\V)\cd\nablag\vphi\ds+\ig\nablag(\nablag\vphi\cd(\V_2+\wu))\cd\nablag\vphi\ds=0.
\end{cases}
\label{tot1}
\end{align}
First, we recall that, by the regularity of strong solutions, there exists $\delta>0$ such that, for $i=1,2$,
$$
\sup_{t\in[0,T]}\norm{\vphi_i}_{L^\infty\gt}\leq 1-\delta,
$$
so that it holds, for $i=1,2$,
\begin{align}
\nonumber&\norm{\V_i}_{H^1_{\L^2}}+\norm{\V_i}_{L^2_{\H^2}}+\norm{\V_i}_{L^\infty_{\H^1_\sigma}}+\norm{\vphi_i}_{L^\infty_{H^3}}+\norm{\vphi_i}_{H^1_{H^1}}\\&+\norm{F''(\vphi_i)}_{L^\infty\gt}+\norm{F'''(\vphi_i)}_{L^\infty\gt}+\norm{F^{(iv)}(\vphi_i)}_{L^\infty\gt}\leq C(T),\label{reguni}
\end{align}
for almost any $t\in[0,T]$. Moreover, concerning the regularity of the flow map, \eqref{regflowmap} holds.
We start form the second identity in \eqref{tot1}. By \eqref{dt3} we have
\begin{align*}
    \frac12\ddt\norm{\nablag\vphi}^2=\ig\na\vphi\cd\na\dt\vphi+\frac12\ig\norma{\na\vphi}^2\Vn 
    H\ds-\ig\na\vphi\cd\E_S(\V_\n)\na\vphi\ds,
\end{align*}
as well as, integrating by parts and recalling that $\Psi''=F''-\theta_0$, 
\begin{align*}
    -\intg \nablag\Delta(\mu_1-\mu_2)\cd\nablag\vphi\ds&=\ig\norma{\na\Delta_\gam\vphi}^2\ds+\ig(F'''(\vphi_1)-F'''(\vphi_2))\norma{\nablag\vphi_1}^2\Delta_{\gam}\vphi\ds\\&+\ig F'''(\vphi_2){\nablag\vphi}\cd\nablag(\vphi_1+\vphi_2)\Delta_\gam\vphi\ds+\ig F''(\vphi_1)\norma{\Delta_\gam\vphi}^2\ds\\&
    +\intg (F''(\vphi_1)-F''(\vphi_2))\Delta_\gam\vphi_2\Delta_\gam\vphi\ds-\intg \theta_0\norma{\Delta_\gam\vphi}^2\ds.
\end{align*}
Furthermore, we have
\begin{align}
    \ig\nablag(\nablag\vphi_1\cd\V)\cd\nablag\vphi\ds=\ig \nablag(\nablag\vphi_1)\V\cd \na\vphi\ds+\ig \nablag^T\V\na\vphi_1\cdot \na\vphi\ds,
\end{align}
and analogously, 
\begin{align}
&\nonumber\ig\nablag(\nablag\vphi\cd(\V_2+\wu))\cd\nablag\vphi\ds\\&=\ig \nablag(\nablag\vphi)(\V_2+\wu)\cd \na\vphi\ds+\ig \nablag^T(\V_2+\wu)\na\vphi\cdot \na\vphi\ds.
\end{align}
Therefore, from \eqref{tot1} we deduce, recalling $F''\geq \theta$,
\begin{align}
    &\nonumber\frac12\ddt\norm{\na\vphi}^2+\norm{\na\Delta_\gam\vphi}^2+\theta\norm{\Delta_\gam\vphi}^2\\&\nonumber\leq \underbrace{\frac12\ig\norma{\na\vphi}^2\vn H\ds-\ig\na\vphi:\E_S(\V_\n)\na\vphi\ds}_{I_1}\nonumber\underbrace{-\ig(F'''(\vphi_1)-F'''(\vphi_2))\norma{\nablag\vphi_1}^2\Delta_{\gam}\vphi\ds}_{I_2}\\&\nonumber\underbrace{-\ig F'''(\vphi_2){\nablag\vphi}\cd\nablag(\vphi_1+\vphi_2)\Delta_\gam\vphi\ds}_{I_3}
    \underbrace{-\intg (F''(\vphi_1)-F''(\vphi_2))\Delta_\gam\vphi_2\Delta_\gam\vphi\ds+\intg \theta_0\norma{\Delta_\gam\vphi}^2\ds}_{I_4}\nonumber\\&
    \underbrace{-\ig \nablag(\nablag\vphi_1)\V\cd \na\vphi-\ig \nablag^T\V\na\vphi_1\cdot \na\vphi\ds}_{I_5}\nonumber\\&
    \underbrace{-\ig \nablag(\nablag\vphi)(\V_2+\wu)\cd \na\vphi-\ig \nablag^T(\V_2+\wu)\na\vphi\cdot \na\vphi\ds}_{I_6}.
    \label{nivers}
\end{align}
We now estimate all these terms separately. First, we have, exploiting \eqref{regflowmap},
\begin{align*}
    \norma{I_1}\leq C(T)\norm{\nablag\vphi}^2.
\end{align*}
Then, exploiting \eqref{reguni} and $H^3\gt\hookrightarrow W^{1,\infty}\gt$, we have, thanks to the strict separation property,
\begin{align*}
\norma{I_2}&= \norma{\ig\int_0^1 F^{(iv)}(s\vphi_1+(1-s)\vphi_2)\vphi ds\norma{\na\vphi_1}^2\Delta_\gam\vphi\ds}\\&\leq C\norm{\vphi_1}^2_{W^{1,\infty}\gt}\norm{\vphi}\norm{\Delta_\gam\vphi}
\leq C(T)\norm{\nablag\vphi}^2+\frac\theta 8 \norm{\Delta_\Gamma\vphi}^2,
\end{align*}
where we used Poincaré's inequality, since $\overline{\vphi}\equiv 0$. Then, similarly,
\begin{align*}
    \norma{I_3}&\leq C\norm{F'''(\vphi_2)}_{L^\infty\gt}\norm{\nablag\vphi}(\norm{\nablag\vphi_1}_{\L^\infty\gt}+\norm{\na\vphi_2}_{\Linfg})\norm{\Delta_\gam\vphi}\\&
    \leq C(T)\norm{\nablag\vphi}^2+\frac\theta 8 \norm{\Delta\vphi}^2.
\end{align*}
Proceeding in the estimate, we have, in a similar way, by the embedding $H^3\gt\hookrightarrow W^{2,4}\gt$,
\begin{align*}
    \norma{I_4}&=\norma{\ig\int_0^1 F'''(s\vphi_1+(1-s)\vphi_2)\vphi ds\Delta_\gam\vphi_2\Delta_\gam\vphi\ds}+\theta_0\norm{\Delta_\gam\vphi}^2\\&\leq C(T)\norm{\vphi}_{L^4\gt}\norm{\Delta_\gam\vphi_2}_{L^4\gt}\norm{\Delta_\gam\vphi}+C\norm{\nablag\Delta_\gam\vphi}\norm{\nablag\vphi}\\&\leq C(T)\norm{\nablag\vphi}^2+\frac\theta8\norm{\Delta_\gam\vphi}^2+\frac12\norm{\na\Delta_\gam\vphi}^2,
\end{align*}
where we used that
\begin{align*}
    \ig \norma{\Delta_\gam\vphi}^2\ds=-\ig \nabla\Delta_\gam\vphi\cdot \nablag\vphi\ds\leq \norm{\nabla\Delta_\gam\vphi}\norm{\na\vphi}.
\end{align*}
About the advective terms, we have, by standard inequalities, \eqref{Korn} and \eqref{Gagliardo}, recalling \eqref{regflowmap} and \eqref{reguni},  
\begin{align*}
\norma{I_5}&\leq   C\norm{\vphi_1}_{W^{2,4}\gt}\norm{\V}_{\L^4\gt}\norm{\nablag\vphi}+C(T)(1+\norm{\V}_{\H^1_\sigma\gt})\norm{\nablag\vphi_1}_{\L^4\gt}\norm{\nablag\vphi}_{\L^4\gt}\\&
\leq 
C(T)(\norm{\V}+\norm{\V}^\frac12\norm{\E_S(\V)}^\frac12)\norm{\nablag\vphi}\\&+C(T)(1+\norm{\V}+\norm{\E_S(\V)})(\norm{\nablag\vphi}+\norm{\nablag\vphi}^\frac12\norm{\Delta_\gam\vphi}^\frac12)\\&
\leq C(T)(\norm{\V}^2+\norm{\na\vphi}^2)+\frac{3}{8}\nu_*\norm{\E_S(\V)}^2+\frac\theta8\norm{\Delta_\gam\vphi}^2.
\end{align*}
Here we have also used the inequality
\begin{align}
\norm{\vphi}_{W^{1,4}\gt}&\leq C(T)\norm{\vphi}_{H^1(\Gamma(t))}^\frac12\norm{\Delta_\gam\vphi}^\frac12\nonumber\\&\leq C(T)\norm{\nablag\vphi}^\frac12(\norm{\vphi}^\frac12+\norm{\Delta_\gam\vphi}^\frac12)\leq C(T)(\norm{\na\vphi}+\norm{\na\vphi}^\frac12\norm{\Delta_\gam\vphi}^\frac12),\label{W14}
\end{align}
which holds by Poincaré's inequality.
In conclusion, by similar estimates, recalling again \eqref{Agmon}, \eqref{regflowmap}, \eqref{reguni}, and \eqref{W14},
\begin{align*}
    \norma{I_6}&\leq C\norm{\vphi}_{H^2\gt}(\norm{\V_2}_{\L^\infty\gt}+1)\norm{\na\vphi}+C(T)(1+\norm{\V_2}_{\H^1_\sigma\gt})\norm{\nablag\vphi_1}_{\L^4\gt}^2\\&
    \leq C(\norm{\vphi}+\norm{\Delta_\gam\vphi})(1+\norm{\V_2}^\frac12\norm{\V_2}_{\H^2\gt}^\frac12)\norm{\na\vphi}+C(T)(\norm{\na\vphi}^2+\norm{\na\vphi}\norm{\Delta_\gam\vphi})\\&
    \leq 
    C(T)(1+\norm{\V_2}_{\H^2\gt}^\frac12+\norm{\V_2}_{\H^2\gt})\norm{\na\vphi}^2+\frac\theta8\norm{\Delta_\gam\vphi}^2.
\end{align*}
Therefore, we have obtained from \eqref{nivers} that  
\begin{align}
&\non\frac12\ddt\norm{\na\vphi}^2+\frac12\norm{\na\Delta_\gam\vphi}^2+\frac{3\theta}{8}\norm{\Delta_\gam\vphi}^2\\&\leq C(T)(1+\norm{\V_2}_{\H^2\gt}^\frac12+\norm{\V_2}_{\H^2\gt})(\norm{\na\vphi}^2+\norm{\V}^2)+\frac{3}{8}\nu_*\norm{\E_S(\V)}^2. 
\label{dtphif}\end{align}
Now we need to analyze \eqref{tot1}$_1$.
First, we notice that 
\begin{align*}
    \frac12\ddt\intg\ro_1\norma{\V}^2\ds =\intg\ro_1\dt\V\cd\V\ds+\frac12\intg \rho'(\vphi_1)\dt\vphi_1\norma{\V}^2\ds+\intg \ro_1\norma{\V}^2\vn H\ds,
\end{align*}
so that, recalling also $\nu\geq \nu_*>0$, we can rewrite \eqref{tot1}$_1$ as
\begin{align}
\non&\frac12\ddt\intg\ro_1\norma{\V}^2\ds+2\nu_*\ig\norma{\E_S(\V)}^2\ds \\&\non\leq \underbrace{\frac12\intg \rho'(\vphi_1)\dt\vphi_1\norma{\V}^2\ds}_{I_7}+\underbrace{\intg \ro_1\norma{\V}^2\vn H\ds}_{I_8}\\&\non-\underbrace{\intg (\rho_1-\rho_2)\dt\V_2\cd\V\ds}_{I_9}-\underbrace{\ig\rho_1(\V\cd\nablag)\V_1\cd\V\ds-\ig(\rho_1-\rho_2)(\V_2\cd\nablag)\V_1\cdot \V\ds}_{I_{10}}\\&\non
    -\underbrace{\ig \ro_2(\V_2\cd\nablag)\V\cd\V\ds}_{I_{11}}-\underbrace{\ig(\J_{\ro_1}\cd\nablag)\V\cd\V\ds-\ig((\J_{\ro_1}-\J_{\ro_2})\cd\nablag\V_2)\cd \V\ds}_{I_{12}}\\&\non
-\underbrace{\ig\ro_1(\wu\cd\nablag)\V\cd \V\ds-\ig \dr(\wu\cd\nablag)\V_2\cd\V\ds-\ig\rho_1(\V\cd\nablag)\wu\cd\V\ds}_{I_{13}}\\&-\non\underbrace{\ig\dr(\V_2\cd\nablag)\wu\cd\V\ds-\ig\ro_1\vn\H\V\cd\V\ds-\ig\dr(\vn\H\V_2-\frac12\nablag(\vn)^2)\cd\V\ds}_{I_{14}}\\&-\non\underbrace{\ig\vn\H(\J_{\ro_1}-\J_{\ro_2})\cd\V\ds}_{I_{15}}-\underbrace{2\ig(\nu(\vphi_1)-\nu(\vphi_2))\E_S(\V_2):\E_S(\V)\ds}_{I_{16}}\\&\non\underbrace{+\ig (\nablag\vphi_1\otimes \nablag\vphi):\nablag\V\ds+\ig (\nablag\vphi\otimes \nablag\vphi_2):\nablag\V\ds}_{I_{17}}\\&\non-\underbrace{2\ig (\nu(\vphi_1)-\nu(\vphi_2))\vn\H:\nablag\V\ds}_{I_{18}}\\&\non
-\underbrace{\ig \dr\dt\wu\cd\V-\ig\dr(\wu\cd\nablag)\wu\cd\V\ds}_{I_{19}}-\underbrace{\ig((\J_{\ro_1}-\J_{\ro_2})\cd\nablag)\wu\cd\V\ds}_{I_{20}}\\&-\underbrace{\ig\dr\vn\H\wu\cd\V\ds-2\ig(\nu(\vphi_1)-\nu(\vphi_2))\E_S(\wu):\E_S(\V)\ds}_{I_{21}}.\label{dtuu} 
\end{align}
We need to estimate all these terms. Now,
by \eqref{Korn} and \eqref{Gagliardo},
\begin{align*}
    \norma{I_7}&\leq C\norm{\dt\vphi_1}\norm{\V}_{\L^4\gt}^2\leq C\norm{\dt\vphi_1}\norm{\V}(\norm{\V}+\norm{\E_S(\V)})
    \leq C\norm{\dt\vphi_1}^2\norm{\V}^2+\nus.
\end{align*}
Then, by \eqref{regflowmap},
\begin{align*}
    \norma{I_8}\leq C(T)\norm{\V}^2.
\end{align*}
Then, recalling that $\rho$ is linear, and again by \eqref{Korn},
\begin{align*}
    \norma{I_9}&\leq C\norm{\vphi}_{L^4\gt}\norm{\dt\V_2}\norm{\V}_{\L^4\gt}\\&\leq C\norm{\na\vphi}\norm{\dt\V_2}\norm{\V}^\frac12(\norm{\V}^\frac12+\norm{\E_S(\V)}^\frac12)\\&
    \leq C(1+\norm{\dt\V_2}^2)(\norm{\na\vphi}^2+\norm{\V}^2)+\nus.
\end{align*}
Then, recalling \eqref{reguni} and Agmon's inequality \eqref{Agmon},
\begin{align*}
    \norma{I_{10}}&\leq C\norm{\V}_{\L^4\gt}\norm{\na\V_1}\norm{\V}_{\L^4\gt}+C\norm{\vphi}_{L^4\gt}\norm{\V_2}_{\Linfg}\norm{\na\V_1}\norm{\V}_{\L^4\gt}\\&
    \leq C(T)\norm{\V}(\norm{\V}+\norm{\E_S(\V)})+C(T)\norm{\na\vphi}\norm{\V_2}^\frac12\norm{\V_2}_{\H^2\gt}^\frac12\norm{\V}^\frac12(\norm{\V}^\frac12+\norm{\E_S(\V)}^\frac12)\\&
    \leq C(T)(1+\norm{\V_2}_{\H^2\gt})(\norm{\V}^2+\norm{\na\vphi}^2)+\nus.
    \end{align*}
    Proceeding in the estimates, we have, similarly,
\begin{align*}
    \norma{I_{11}}&\leq C\norm{\V_2}_{\Linfg}\norm{\nablag\V}\norm{\V}\leq C\norm{\V_2}^\frac12\norm{\V_2}_{\H^2\gt}^\frac12(\norm{\V}+\norm{\E_S(\V)})\norm{\V}\\&
    \leq C(\norm{\V_2}_{\H^2\gt}^\frac12+\norm{\V_2}_{\H^2\gt})\norm{\V}^2+\nus.
\end{align*}
Now, recalling the definition of $\J_{\rho_i}$, we have 
\begin{align}
    \J_{\rho_1}-\J_{\rho_2}&=-\frac{\trho_1-\trho_2}{2}\nablag(\mu_1-\mu_2)\nonumber\\&=\frac{\trho_1-\trho_2}{2}\nablag\Delta_\gam\vphi-\frac{\trho_1-\trho_2}{2}(F''(\vphi_1)-\theta_0)\na\vphi-\frac{\trho_1-\trho_2}{2}(F''(\vphi_1)-F''(\vphi_2))\na\vphi_2,\label{jro}
\end{align}
so that, since $F''(\vphi_1)-F''(\vphi_2)=\int_0^1 F'''(s\vphi_1+(1-s)\vphi_2)\vphi ds$, recalling \eqref{Korn}, \eqref{Gagliardo}, Agmon's inequality \eqref{Agmon}, \eqref{reguni} and the strict separation property,
\begin{align*}
   & \norma{I_{12}}\leq C\norm{\nablag\mu_1}_{\L^\infty\gt}\norm{\na\V}\norm{\V}\\&+C\norm{\na\Delta_\gam\vphi}\norm{\nablag\V_2}_{\L^4\gt}\norm{\V}_{\L^4\gt}+C\norm{F''(\vphi_1)-\theta_0}_{L^\infty\gt}\norm{\na\vphi}\norm{\nablag\V_2}_{\L^4\gt}\norm{\V}_{\L^4\gt}\\&
   +C\norm{\na\vphi_2}_{\L^\infty\gt}\norm{\vphi}\norm{\nablag\V_2}_{\L^4\gt}\norm{\V}_{\L^4\gt}\\&
   +C\norm{\na\mu_1}^\frac12\norm{\mu_1}_{H^3\gt}^\frac12(\norm{\V}+\norm{\E_S(\V)})\norm{\V}\\&+C\norm{\na\Delta_\gam\vphi}\norm{\nablag\V_2}^\frac12\norm{\V_2}_{\H^2\gt}^\frac12\norm{\V}^\frac12(\norm{\V}^\frac12+\norm{\E_S(\V)}^\frac12)\\&
   +C\norm{\nablag\vphi}\norm{\nablag\V_2}^\frac12\norm{\V_2}_{\H^2\gt}^\frac12\norm{\V}^\frac12(\norm{\V}^\frac12+\norm{\E_S(\V)}^\frac12)\\&
   +C\norm{\vphi}\norm{\nablag\V_2}^\frac12\norm{\V_2}_{\H^2\gt}^\frac12\norm{\V}^\frac12(\norm{\V}^\frac12+\norm{\E_S(\V)}^\frac12)\\&
   \leq 
   C(T)(1+\norm{\mu_1}_{H^3\gt}+\norm{\V_2}_{\H^2\gt}^\frac12+\norm{\V_2}_{\H^2\gt})(\norm{\V}^2+\norm{\nablag\vphi}^2)\\&+\frac{1}{4}\norm{\na\Delta_\gam\vphi}^2+\nus.
\end{align*}
Furthermore, recalling \eqref{regflowmap},
\begin{align*}
    \norma{I_{13}}&\leq C(T)(\norm{\wu}_{\Linfg}\norm{\nablag\V}\norm{\V}+\norm{\vphi}_{L^4\gt}\norm{\wu}_{\Linfg}\norm{\na\V_2}\norm{\V}_{\L^4\gt}+\norm{\V}^2\norm{\wu}_{\Linfg})
    \\&\leq C(T)(\norm{\V}^2+\norm{\na\vphi}^2)+\nus.
\end{align*}
Then, similarly 
\begin{align*}
    &\norma{I_{14}}\\&\leq C\norm{\vphi}_{L^4   \gt}\norm{\V_2}_{\Linfg}\norm{\na\wu}_{\Linfg}\norm{\V}_{\L^4\gt}+C(T)\norm{\V}^2+C(T)\norm{\vphi}(1+\norm{\V_2}_{\Linfg})\norm{\V} 
    \\&
    \leq C\norm{\na\vphi}\norm{\V_2}_{\H^2\gt}^\frac12\norm{\V}^\frac12(\norm{\V}^\frac12+\norm{\E_S(\V)}^\frac12)+C(T)\norm{\V}^2+C(T)\norm{\na\vphi}\norm{\V_2}_{\H^2\gt}^\frac12\norm{\V}\\&
    \leq C(T)(1+\norm{\V_2}_{\H^2\gt}^\frac12+\norm{\V_2}_{\H^2\gt})(\norm{\V}^2+\norm{\na\vphi}^2)+\nus,
\end{align*}
as well as, recalling \eqref{reguni} and \eqref{jro},
\begin{align*}
    &\norma{I_{15}}+\norma{I_{20}}\leq C(\norm{\vn\H}_{\Linfg}+\norm{\na\wu}_{\Linfg})\norm{\J_{\rho_1}-\J_{\rho_2}}\norm{\V}\\&
    \leq C(T)(\norm{\na\Delta_\gam\vphi}+\norm{F''(\vphi)-\theta}_{L^\infty\gt}\norm{\na\vphi}+\norm{\vphi}\norm{\nablag\vphi_2}_{\Linfg})\norm{\V}\\&
    \leq C(T)(\norm{\V}^2+\norm{\na\vphi}^2)+\frac18 \norm{\na\Delta_\gam\vphi}^2.
\end{align*}
Then, thanks to the regularity of $\nu$ and again to \eqref{Gagliardo}, \eqref{regflowmap} and \eqref{reguni},
\begin{align*}
\norma{I_{16}}+\norma{I_{21}}&\leq 
C\norm{\vphi}_{L^4\gt}(\norm{\V_2}_{\W^{1,4}\gt}+\norm{\wu}_{\W^{1,4}\gt})\norm{\E_S(\V)}\\&\leq 
C(T)\norm{\na\vphi}(1+\norm{\V_2}^\frac12_{\H^1\gt}\norm{\V_2}^\frac12_{\H^2\gt})\norm{\E_S(\V)}\\&
\leq C(T)(1+\norm{\V_2}_{\H^2\gt})\norm{\na\vphi}^2+\nus.
\end{align*}
Furthermore, concerning $I_{17}$, recalling \eqref{Korn}, \eqref{Gagliardo}, and \eqref{reguni},
\begin{align*}
    \norma{I_{17}}&\leq \norm{\na\vphi}_{\L^4\gt}(\norm{\na\vphi_1}_{\L^4\gt}+\norm{\na\vphi_2}_{\L^4\gt})(\norm{\V}+\norm{\E_S(\V)}\\&
    \leq C(T)\norm{\na\vphi}^\frac12(\norm{\na\vphi}^\frac12+\norm{\Delta_\gam\vphi}^\frac12)(\norm{\V}+\norm{\E_S(\V)}\\&
    \leq C(T)(\norm{\V}^2+\norm{\na\vphi}^2)+\nus+\frac{\theta}{8}\norm{\Delta_\gam\vphi}^2.
\end{align*}
To conclude, terms $I_{18}$ and $I_{19}$ are straightforwardly estimated thanks to Poincaré's inequality and \eqref{regflowmap}, since
\begin{align*}
&\norma{I_{18}}+\norma{I_{19}}\\&\leq C\norm{\vphi}\norm{v_\n\H}_{\Linfg}\norm{\na\V}+C\norm{\vphi}\norm{\dt\wu}_{\Linfg}\norm{\V}+C\norm{\vphi}\norm{\wu}_{\Linfg}\norm{\na\wu}_{\Linfg}\norm{\V}\\&
\leq C(T)(\norm{\na\vphi}^2+\norm{\V}^2)+\nus.
\end{align*}
Putting all these estimates together, we deduce from \eqref{dtuu} that
\begin{align*}
&\frac12\ddt\intg\ro_1\norma{\V}^2\ds+\frac{11}{8}\nu_*\ig\norma{\E_S(\V)}^2\ds\\&\leq C(T)(1+\norm{\dt\vphi_1}^2+\norm{\dt\V_2}^2+\norm{\V}_{\H^2\gt}+\norm{\mu_1}_{H^3\gt})(\norm{\V}^2+\norm{\na\vphi}^2)\\&+\frac{3}{8}\norm{\nablag\Delta_\gam\vphi}^2+\frac\theta 8\norm{\Delta_\gam\vphi}^2,
\end{align*}
which, if we sum up with \eqref{dtphif} and recall that $\rho_1\geq \rho_*$, becomes
\begin{align*}
&\frac12\ddt\left(\intg\left(\ro_1\norma{\V}^2+\norm{\na\vphi}^2\right)\ds\right)+\nu_*\ig\norma{\E_S(\V)}^2\ds+\frac18\norm{\na\Delta_\gam\vphi}^2+\frac{\theta}4\norm{\Delta_\gam\vphi}^2\\&\leq C(T)(1+\norm{\dt\vphi_1}^2+\norm{\dt\V_2}^2+\norm{\V_2}_{\H^2\gt}+\norm{\mu_1}_{H^3\gt})\left(\ig\ro_1\norma{\V}^2\ds+\norm{\na\vphi}^2\right).
\end{align*}
Now, thanks to \eqref{reguni}, it holds 
$$
\mathcal{H}:=1+\norm{\dt\vphi_1}^2+\norm{\dt\V_2}^2+\norm{\V_2}_{\H^2\gt}+\norm{\mu_1}_{H^3\gt}\in L^1(0,T),
$$
entailing by the Gronwall Lemma the desired continuous dependence estimate, leading to the uniqueness of strong solutions. This concludes the proof of the theorem.

\medskip
\textbf{Acknowledgments.} Part of this work was done while AP was visiting HA and HG at the Department of Mathematics of the University of Regensburg, whose hospitality is kindly acknowledged. The authors gratefully acknowledge the support by the Graduiertenkolleg 2339 IntComSin of the Deutsche
Forschungsgemeinschaft (DFG, German Research Foundation) – Project-ID 321821685. 
AP is a member of Gruppo Nazionale per l’Analisi Matematica, la Probabilità e le loro Applicazioni (GNAMPA) of
Istituto Nazionale per l’Alta Matematica (INdAM).  This research was funded in part by the Austrian Science Fund (FWF) \href{https://doi.org/10.55776/ESP552}{10.55776/ESP552}.
For open access purposes, the authors have applied a CC BY public copyright license to
any author accepted manuscript version arising from this submission.

\bibliography{Bibliography}
\bibliographystyle{abbrv}
\end{document}